\documentclass{article}    
\newcommand{\QED}{{\par\hfill$\square$\par}}

\newcommand{\Div}{\mbox{\rm div}\,}

\newcommand{\curl}{\mbox{\rm curl}\,}

\newcommand{\Int}[2]{{\displaystyle \int_{ #1}^{ #2}}}

\newcommand{\Frac}[2]{\displaystyle{\frac{\displaystyle{#1}}{\displaystyle{#2}}}}

\newcommand{\beea}{\begin{eqnarray}}

\newcommand{\eeea}{\end{eqnarray}}

\newcommand{\bfe}{{\mbox{\boldmath $e$}} }

\newcommand{\0}{{\mbox{\boldmath $0$}} }

\newcommand{\Bd}{\begin{defn}\begin{rm}}
\newcommand{\EED}[1]{\end{rm}\label{defn:#1}\end{defn}}

\newcommand{\BF}{\begin{footnotesize}}

\newcommand{\EF}{\end{footnotesize}}

\newcommand{\ode}[2]{{\displaystyle \frac{\mbox{$d #1$}}{\mbox{$d #2$}}}}

\newcommand{\bi}{\begin{itemize}}

\newcommand{\ei}{\end{itemize}}

\newcommand{\ed}{\end{document}}

\newcommand{\be}{\begin{equation}}

\newcommand{\ba}{\begin{array}}

\newcommand{\ea}{\end{array}}

\newcommand{\ee}{\end{equation}}

\newcommand{\eeq}[1]{\label{eq:#1}\end{equation}}

\newcommand{\real}{{\mathbb R}}

\newcommand{\bfomega}{\mbox{\boldmath $\omega$}}

\newcommand{\bfx}{\mbox{\boldmath $x$}}

\newcommand{\bfv}{{\mbox{\boldmath $v$}} }

\newcommand{\bfu}{{\mbox{\boldmath $u$}} }

\newcommand{\bfw}{{\mbox{\boldmath $w$}} }

\newcommand{\bfa}{{\mbox{\boldmath $a$}} }

\newcommand{\bfA}{{\mbox{\boldmath $A$}} }

\newcommand{\bfM}{{\mbox{\boldmath $M$}} }

\newcommand{\bfB}{{\mbox{\boldmath $B$}} }

\newcommand{\bfG}{{\mbox{\boldmath $G$}} }

\newcommand{\bfI}{{\mbox{\boldmath $I$}} }

\newcommand{\bfN}{{\mbox{\boldmath $N$}} }

\newcommand{\bfh}{{\mbox{\boldmath $h$}} }

\newcommand{\calc}{{\cal C}}

\newcommand{\calf}{{\cal F}}

\newcommand{\calp}{{\cal P}}

\newcommand{\calq}{{\cal Q}}

\newcommand{\cals}{{\cal S}}

\newcommand{\bfT}{{\mbox{\boldmath $T$}} }

\newcommand{\bfF}{{\mbox{\boldmath $F$}} }

\newcommand{\bfL}{{\mbox{\boldmath $L$}} }

\newcommand{\bfg}{{\mbox{\boldmath $g$}} }

\newcommand{\bfn}{{\mbox{\boldmath $n$}} }

\newcommand{\half}{\mbox{$\frac{1}{2}$}}

\def\Bbb R{\real}

\def\hat{\widehat}

\def\tilde{\widetilde}

\def\bar{\overline}

\newcommand{\bfgamma}{\mbox{\boldmath $\gamma$}}

\newcommand{\ED}{\end{description}}

\newcommand{\defref}[1]{{\rm Definition \ref{defn:#1}}}

\newcommand{\Br}{\begin{rem}\begin{rm}}

\newcommand{\Er}{\end{rm}\end{remark}}

\newtheorem{lemm}{Lemma}[section]

\newtheorem{theo}{Theorem}[section]

\newtheorem{prop}{Proposition}[section]

\newtheorem{rem}{Remark}[section]

\newtheorem{defn}{Definition}[section]

\newtheorem{coro}{Corollary}[section]

\newtheorem{exe}{\footnotesize{Exercise}}[section]

\newcommand{\Be}{\begin{exe}\begin{footnotesize}\begin{rm}}

\newcommand{\EE}[1]{\end{rm}\end{footnotesize}\label{exe:#1}\end{exe}}

\newcommand{\Bt}{\begin{theo}\begin{sl}}

\newcommand{\Bp}{\begin{prop}\begin{sl}}

\newcommand{\EP}[1]{\end{sl}\label{prop:#1}\end{prop}}

\newcommand{\Et}{\end{sl}\end{theorem}}

\newcommand{\Bl}{\begin{lemm}\begin{sl}}

\newcommand{\El}{\end{sl}\end{lemma}}

\newcommand{\eqref}[1]{{\rm (\ref{eq:#1})}}

\newcommand{\Bc}{\begin{coro}\begin{sl}}

\newcommand{\Ec}{\end{sl}\end{coro}}

\newcommand{\ET}[1]{\end{sl}\label{theo:#1}\end{theo}}

\newcommand{\EL}[1]{\end{sl}\label{lemm:#1}\end{lemm}}

\newcommand{\theoref}[1]{{\rm Theorem \ref{theo:#1}}}

\newcommand{\propref}[1]{{\rm Proposition \ref{prop:#1}}}

\newcommand{\ER}[1]{\end{rm}\label{rem:#1}\end{rem}}

\newcommand{\EC}[1]{\end{sl}\label{coro:#1}\end{coro}}

\newcommand{\remref}[1]{{\rm Remark \ref{rem:#1}}}

\newcommand{\lemmref}[1]{{\rm Lemma \ref{lemm:#1}}}

\usepackage{amssymb}
\usepackage{amsfonts}

\renewcommand{\real}{{\mathbb R}}

\newcommand{\boms}{{\bfomega_*}}
\usepackage[mathscr]{eucal}
\usepackage{color}
\usepackage{sectsty}
\allsectionsfont{\bfseries\sffamily}
\begin{document}
\title{Stability of Permanent Rotations and Long-Time Behavior of  Inertial Motions of a Rigid Body with an Interior Liquid-Filled Cavity} 
\author{Giovanni P. Galdi}
\date{University of Pittsburgh, USA}
\maketitle
\pagenumbering{arabic}

\begin{abstract}A rigid body, with an interior cavity entirely filled with a Navier-Stokes liquid,  moves in absence of external torques relative to the center of mass, $G$, of the coupled system body-liquid (inertial motions). The only steady-state motions allowed about $G$ are then those where the system, as a whole rigid body, rotates uniformly around one of the central axes of inertia (permanent rotations). Objective of this article is two-fold. On the one hand, we provide sufficient conditions for the asymptotic, exponential stability of permanent rotations, as well as for their instability. On the other hand, we study the asymptotic behavior of the generic motion in the class of weak solutions and show that there exists a time $t_0$ after that all such solutions must decay exponentially fast to a permanent rotation. This result provides a {\em full} and rigorous explanation of Zhukovsky's conjecture, and explains, likewise, other interesting phenomena that are observed in both lab and numerical experiments.                                                 

\end{abstract}
             
\section*{\normalsize Introduction}
The problem of the motion of the coupled system constituted by a rigid body with an interior cavity that is entirely filled with a liquid has represented, over the centuries, one of the main focuses of theoretical and applied research. As a matter of fact, the first mathematical analysis of such a problem can be traced  back to the pioneering work of Stokes concerning the motion of a rectangular box filled  with an inviscid liquid \cite[\S 13]{Stok}. 
\par
One of the main reasons why  this topic received  all along increasing attention is because it was rather immediately recognized that the dynamics of the rigid body can be substantially and drastically affected by the presence of the liquid, in several different and major aspects.
In these regards, the  finding of { Lord Kelvin} ({W. Thomson}) can be considered a true cornerstone. Actually, 
it was well known from both a theoretical viewpoint  --basically, by Lagrange \cite[Part 2, \S 9]{La} and Poisson \cite{Poi}-- and practical one --after the use of gyroscopes in navigation since the early 1740s-- that  uniform rotations occurring around either the shorter  or the longer  axis of a spheroid, in absence of external torques ({\em inertial motions}), are  {\em both stable}.
However, Kelvin's experiment showed that if a  thin--walled
spheroidal gyroscope  is filled  up with water ({\em ``liquid gyrostat"}), it  would be stable when set in rotation around its shorter axis, whereas it would be unstable otherwise, no matter how large the magnitude of the initial angular velocity. Notice that the shorter axis is the one with respect to which the moment of inertia of the system is a maximum
The mathematical explanation of Kelvin's experiment was the object of studies by several prominent mathematicians, including Poincar\'e \cite{Po} --who  considered the more general case of an elastic container-- and Basset [13], who analyzed the
homogeneous vortex motion of a liquid in an ellipsoidal cavity. However, the outcome of these investigations provided only  approximate or incomplete answers. 
\par
Another fundamental insight to the problem is due to N.Y. Zhukovsky. In his thorough analysis  \cite{Zhu}, Zhukovsky puts forward a completely unexpected property of the coupled system body-liquid, $\mathscr S$, in the case when the liquid filling the cavity is {\em viscous} and $\mathscr S$ moves by inertial motion. More precisely, on the basis of a straightforward energy analysis, he envisages that  the liquid should produce a substantial {\em stabilizing effect} on the motion of the body in a way that the terminal state should be one where $\mathscr S$ moves as a single rigid body by  uniform, rotational motion around one of the central axes of inertia \cite[\S 38]{Zhu}. 
This dynamical behavior is entirely at odds with the one that the body might perform with an {\em empty} cavity, where  the unsteady motion is much more complicated and, of course, far from reaching any steady-state configuration (e.g., \cite[\S 4]{Lei}).  It must be emphasized, at this point, that Zhukovsky's argument is altogether of heuristic nature and, therefore,  lacks of sound mathematical rigor. As a result, even though the above property  is often referred to (especially in the Russian literature, e.g., \cite[p. 98]{MR}, \cite[p. 3]{Ch}) as ``Zhukovsky's theorem",  it is more precise to call it, instead, ``Zhukovsky's conjecture." However, \textcolor{black}{despite} the absence of a rigorous analysis,  the use of interior cavities filled with liquid as dampers in rigid and elastic  structures is a common  procedure,  adopted since the mid 1960s \cite{Co} in different branches of applied sciences, especially space technology (e.g., \cite{Bu}) and civil engineering (e.g., \cite{Aly}).    

Coming to the mathematical analysis and interpretation  of the phenomena described above and, more generally,  the generic motion of a rigid body with a liquid-filled cavity, the classical literature includes a very large number of contributions. However, this body of work, probably also due to the strong influence of the seminal articles of Rumyantsev \cite{Ru} and Chernusko \cite{Ch} on the entire field, is only seldom of rigorous nature, and mostly based  on a simplified set of equations --that at times reduces to ordinary differential equations-- and/or special shapes of the body and cavity. Since it would be hopeless to cite all the relevant literature, we refer the reader to the monographs \cite{MR,Ch,kk,kk1} and the bibliography there cited.       

Very recently, the present author, jointly with his associates, has started a systematic study of the motion of a rigid body with a liquid-filled cavity \cite{M,GMZ0,GMZ,DGMZ,GaMaMa,Ma2}. 
In these papers, one  main  objective, among others,  was to provide a mathematically rigorous explanation of Kelvin's experiment and a likewise solid proof of Zhukovsky's conjecture.
The outcome of this effort has been remarkably successful \cite{DGMZ} though, however, not entirely satisfactory. In fact, on the one hand, Kelvin's experiment is only partially recovered. Actually, just in the case of systems such as prolate spheroids entirely filled with a liquid, in \cite[Theorem 6(e)]{DGMZ} it is  shown, in particular, that permanent rotations around the longer axis are  {\em unstable},  but not that those around the shorter axes are indeed {\em stable}, as demonstrated by  Kelvin's finding. On the other hand, in \cite[Theorem 4]{DGMZ} Zhukovskiy's conjecture is rigorously proved for all types of body-liquid systems, with the exception of those whose mass distribution is such that two central moments of inertia coincide and are strictly greater than the third one. This happens, for example, in cylindrically-shaped containers filled with liquid, when the radius of the base is shorter than the height (like in a soda can). 
Finally,  the results presented in \cite{DGMZ} are not able to explain  another interesting phenomenon that both lab  \cite{WJ} and numerical \cite[Section 9.1]{DGMZ} experiments strongly suggest. More precisely, it is observed that after a finite interval of time, whose length mostly depends on the viscosity of the liquid and the size of the initial conditions, the coupled system almost abruptly reaches a uniform terminal state. In fact, both experiments indicate that once the viscosity of the liquid has reduced its motion ``sufficiently close" to the  (relative) rest, the rate of decay to the terminal state appears  be of exponential type.        
\par
The main goal of this article is to analyze in  details the stability properties of {uniform rotations} ({\em ``permanent rotations"}) and long-time behavior of  motions of a rigid body with an interior cavity entirely filled with a viscous liquid around its center of mass $G$, in absence of external torques relative to $G$ \textcolor{black}{({\em ``inertial motions"})}. In doing so, we shall, in particular, provide a positive answer to all problems left open in \cite{DGMZ} and mentioned  above. Besides  $C^2$-smoothness of the cavity, we do not make any other assumption about its shape or the shape of the body. \par Our approach is quite different than the one adopted in \cite{DGMZ}, which is based on (appropriately modified) tools borrowed from classical dynamical system theory. Instead, the method we use here relies upon a detailed study of the spectrum of the relevant linear operator, $\bfL$, obtained by linearizing the full nonlinear operator around a given permanent rotation, ${\sf s}_0$, of $\mathscr S$. As is well known, this rotation may only occur about  an axis, $\bfe$, coinciding with one of  the eigenvectors,  $\bfe_1,\bfe_2$, and $\bfe_3$, of the inertia tensor    of $\mathscr S$ relative to $G$ ({\em central axes of inertia}). Let  $A,B$, and $C$ denote, in the order, the associated eigenvalues  ({\em central moments of inertia}), and for $\lambda\in\{A,B,C\}$, by $\cals(\lambda)$ the corresponding eigenspace.
Without loss of generality, we take $A\le B\le C$.  We then show (\lemmref{1.2}) that, for any given ${\sf s}_0$, the spectrum $\sigma(\bfL)$ of $\bfL$ is purely discrete with  eigenvalues clustering only at infinity. Furthermore, 0 is always an eigenvalue with algebraic multiplicity $m={\rm dim}\,(\cals(\lambda))\in \{1,2,3\}$.  This implies, in particular, the existence of a center manifold, $\mathscr C$, in our problem. However, we show (\propref{1.2}) that $\mathscr C$ is ``slow'', namely, 0 is the only point of $\sigma(\bfL)$ on the imaginary axis. In addition, we prove that  $0$ is \textcolor{black}{semisimple} and this allows us to characterize the sign of $\Re[\sigma(\bfL)]\backslash \{0\}$ in terms of the central moments $A,B$, and $C$ and the direction $\bfe$ (\propref{1.2}). We then employ this information along with a form of the ``generalized linearization principle"  \cite{PSZ} to show that  if the uniform rotation ${\sf s}_0$ occurs around an axis with maximum moment of inertia, then it is asymptotically {\em exponentially} stable; see \propref{1.1} and \theoref{1.1}. By this we mean that  ${\sf s}_0$  is stable in the sense of Lyapunov (in appropriate norms), and, moreover, every motion starting in a suitable neighborhood of ${\sf s}_0$ will converge exponentially fast to a terminal state that is still a uniform rotation around $\bfe$; see \defref{1.1}. The latter, however, will in general  be different from ${\sf s}_0$, due to the conservation of total angular momentum for $\mathscr S$; see \remref{3.1}. Conversely, we show that if ${\sf s}_0$ occurs around an axis of minimum moment of inertia, then it is unstable (\propref{1.1} and \theoref{1.1}), and so furnishing, in particular, a full explanation of the outcome of Kelvin's experiment; see \remref{3.2}.  \par 
The spectral properties of $\bfL$ mentioned earlier on also enable us  to give a rigorous proof of the ``abrupt" decay of the coupled system body-liquid to the terminal state of uniform rotation. In fact, also with the help of the results already established in \cite{DGMZ}, we prove that there is a time $t_0> 0$ after that all solutions possessing  finite kinetic energy at time $t=0$ ({\em weak solutions}) must decay {\em exponentially fast} to their terminal state. This result provides a full proof of Zhukovsky's conjecture; see \theoref{2.1}. Finally, we investigate the question of around which central axis the terminal uniform rotation will take place ({\em attainability problem}). Again with the help of the results established in \cite{DGMZ} we show that, for an open set of initial data, all corresponding weak solution after some time $t_0>0$ will converge at an {\em exponential rate} to a rotation around the axis with respect to which the moment of inertia is a maximum; see \theoref{5.1}.  The norm in which this convergence is established is quite strong, since it involves first time derivatives and second spatial derivatives of the velocity field of the liquid, and angular velocity and acceleration of the body.  \par
The plan of the paper is as follows. In Section 1 we give the mathematical formulation of  inertial motions of a body with an interior cavity entirely filled with a Navier-Stokes liquid. We then show  that the associated perturbation problem can be formulated as an abstract evolution equation in a suitable Hilbert space and show, by using a classical semigroup approach, that the initial-value problem possesses a unique smooth solution, at least locally in time. Section 2 is devoted to the study of the spectrum of the relevant linear operator $\bfL$, for which we show the properties reported earlier on. These results are then employed in the  following Section 3 to show, on the one hand, that the local solution constructed in Section 1 can be made global in time for sufficiently ``small" data, and, on the other hand, to provide  sufficient conditions for  the exponential stability and instability of permanent rotations. In Section 4 we investigate the long-time behavior of the generic motion and show that, under suitable assumptions on the central moments of inertia, any weak solution after some time $t_0$ must converge exponentially fast to a permanent rotation (Zhukowsky's conjecture). A crucial role in the proof is played by \lemmref{2.1} that ensures that the velocity field of {\em any} weak solution must decay to 0, as times diverges, in higher-order spatial norms. The final Section 5 is dedicated to the problem of attainability of permanent rotations.      
\setcounter{section}{1}
\renewcommand{\theequation}{{1}.\arabic{equation}}
\section*{\normalsize 1. Formulation of the Problem and Local Existence Theory}
The motion of the coupled system body-liquid, $\mathscr S$, is governed by the following set of equations (e.g., \cite{DGMZ}):
\be\ba{cc}\medskip
\left.\ba{ll}\medskip
\bfv_t+\bfv\cdot\nabla\bfv+\dot{\bfomega}\times\bfx+2\bfomega\times\bfv=\nu\Delta\bfv-\nabla p\\
\Div\bfv=0\ea\right\}\ \ \mbox{in $\calc\times (0,\infty)$}\\

\bfv(x,t)=\0\ \ \mbox{at $\partial\calc$}
\ea\eeq{1.1}
and
\be
\dot{\bfM}+\bfomega\times\bfM=\0\,,\ \ \bfM:=\mathbb I\cdot\bfomega+\int_\calc \bfx\times\bfv\,.
\eeq{1.2}
Here, $\bfv$ and $\rho p$ (with $\rho$ constant density of the liquid that, without loss of generality, will be taken to be 1 for simplicity) are velocity and (modified) pressure field of the liquid,  $\bfomega$ is the angular velocity of the body, and  $\calc$ is the domain of $\real^3$ occupied by the liquid. Throughout this paper we shall assume that $\calc$ is of class $C^2$. Moreover, $\mathbb I$ denotes the inertia tensor of the coupled system $\mathscr S$. We finally notice that the vector function $\bfM$ represents the total angular momentum of $\mathscr S$. 

Equation \eqref{1.1}--\eqref{1.2} are written in a body-fixed frame, $\calf$ with the origin at the center of mass, $G$, of $\mathscr S$. We choose
$\calf\equiv\{G,\bfe_i\}$ where $\{\bfe_i\}$, are eigenvectors of the inertia tensor $\mathbb I$ (central axes of inertia). Moreover, we denote by $A,B$, and $C$ the central moments of inertia, namely, the eigenvalues of $\mathbb I$ corresponding to the eigenvectors  $\bfe_1,\bfe_2$, and $\bfe_3$, respectively. Without loss of generality, we shall assume $A\le B\le C$.  
\par
For $\lambda\in\{A,B,C\}$ we let $\cals(\lambda)$ be the eigenspace associated to $\lambda$. If we pick a unit vector $\bfe\in\cals(\lambda)$ it is immediately verified that
\be
{\sf s}_0=(\bfv_0\equiv \0,\bfomega_0=\omega_0\bfe)\,,\ \ \omega_0\in\real-\{0\}\,,
\eeq{1.0}
is a (non-trivial) steady-state solution to  \eqref{1.1}--\eqref{1.2} representing a {\em permanent rotation} performed by the coupled system as a single rigid body. In fact, permanent rotations are the only steady-state motions allowed for the coupled system $\mathscr S$,  and, as is well known, can only occur around a central axis of inertia.
\par 
Even though of rather trivial proof, for the relevance acquired later on, we would like to single out the following result in the form of a lemma.  
\Bl Suppose $\bfomega_0\in \cals(\lambda)-\{\0\}
$. Then, $\bfomega\in\real^3$ is a solution to
\be
\bfomega_0\times\mathbb I\cdot\bfomega+\bfomega\times\mathbb I\cdot\bfomega_0=\0\,.
\eeq{00}
if and only if $\bfomega\in\cals(\lambda)$.
\EL{1.0}
{\em Proof.} By assumption,
\be\mathbb I\cdot\bfomega_0=\lambda\,\bfomega_0\,.\eeq{001}
Now, if $\bfomega\in \cals(\lambda)$, we also have
$$
\mathbb I\cdot\bfomega=\lambda\,\bfomega\,.
$$
so that, by \eqref{001},
$$
\bfomega_0\times\mathbb I\cdot\bfomega+\bfomega\times\mathbb I\cdot\bfomega_0=\lambda\,(\bfomega_0\times\bfomega+\bfomega\times\bfomega_0)=\0\,.
$$
Conversely,  from \eqref{00} and \eqref{001} we get
$$
\mathbb I\cdot\bfomega-\lambda\,\bfomega=\sigma\,\bfomega_0\,\ \ \mbox{some $\sigma\in\real$}\,.
$$
Dot-multiplying both sides of this equation by $\bfomega_0$ and using  the symmetry of $\mathbb I$ along with \eqref{001}, we deduce $\sigma\,|\bfomega_0|^2=0$,
which is possible if and only if $\sigma=0$.
\QED

One of our main objectives is to determine necessary and sufficient conditions for the stability of the steady-state solution \eqref{1.0}, in a sense made precise later on in \defref{1.1}. To this end, 
denoting by $(\bfv,p,\bfomega_0+\bfomega)$ a generic solution to  \eqref{1.1}--\eqref{1.2} with $\bfomega_0\in\cals(\lambda)$, we  deduce that the ``perturbation" $(\bfv,p,\bfomega)$ must satisfy the following system
\be\ba{cc}\medskip
\left.\ba{ll}\medskip
\bfv_t+\dot{\bfomega}\times\bfx+2\bfomega_0\times\bfv-\nu\Delta\bfv-\nabla p=- 2\bfomega\times\bfv-\bfv\cdot\nabla\bfv\\
\Div\bfv=0\ea\right\}\ \ \mbox{in $\calc\times(0,\infty)$}\\

\bfv(x,t)=\0\ \ \mbox{at $\partial\Omega$}
\ea\eeq{1.3}
and
\be
\mathbb I\cdot\dot{\bfomega}-\mathbb I\cdot\dot{\bfa}+\bfomega_0\times\mathbb I\cdot\bfomega+\bfomega\times\mathbb I\cdot\bfomega_0-\bfomega_0\times\mathbb I\cdot\bfa=-\bfomega\times\mathbb I\cdot\bfomega+\bfomega\times\mathbb I\cdot\bfa\,,
\eeq{1.4}
where
\be
\bfa:=-\mathbb I^{-1}\cdot\int_\calc \bfx\times\bfv\,.
\eeq{1.5}
We shall now rewrite \eqref{1.3}--\eqref{1.4} in a suitable abstract form. To this end, 
we introduce the Hilbert space\footnote{The notation used in this article is quite standard; see e.g. \cite{IVP}.}  
$$
H:=L^2_\sigma(\calc)\oplus \real^3=\big\{\bfu=(\bfv,\bfomega)^\top: \bfv\in L^2_\sigma(\calc),  \ \bfomega\in \real^3\big\}\,,
$$
where
$$
L^2_\sigma(\calc)=\{\bfv\in L^2(\Omega):\ \Div\bfv=0\ \mbox{in $\calc$\,, and $\bfv\cdot\bfn|_{\partial\calc}=0$}\}
\,,
$$
and $\bfn$ is the unit exterior normal on $\partial\calc$. 
The scalar product of two elements  $\bfu_i=(\bfv_i,\bfomega_i)^\top$, $i=1,2$, in $H$ is defined by
$$
\langle \bfu_1,\bfu_2\rangle:=\int_\calc\bfv_1\cdot\bfv_2\,{ d}\calc+\bfomega_1\cdot\bfomega_2\,,\ \
$$
with associated
norm
$$
\|\bfu\|:=\langle\bfu,\bfu\rangle^{\frac12}\,.
$$
We next define the following operators
\be\ba{ll}\medskip
\bfI: \bfu\in H\mapsto \bfI\bfu:=\Big( \bfv+{\rm P}\,(\bfomega\times \bfx),\ \mathbb I\cdot(\bfomega-\bfa)\Big)^\top\in H\,,\\ \medskip
\bfA:\bfu\in D(\bfA):=[W^{2,2}(\calc)\cap W_0^{1,2}(\calc)\cap L^2_\sigma(\calc)]\oplus \real^3\subset H\\ \medskip
\hspace*{7cm}\mapsto \bfA\bfu:=\big(-\nu\,{\rm P}\,\Delta\bfv\,,\ \bfomega)^\top\in H\,,\\ \medskip
\bfB: \bfu\in H\mapsto \bfB\bfu:=\Big(2\,{\rm P}\,(\bfomega_0\times\bfv)\,,\ \bfomega_0\times\mathbb I\cdot\bfomega+\bfomega\times\mathbb I\cdot\bfomega_0-\bfomega_0\times\mathbb I\cdot\bfa-\bfomega\Big)^\top\in H\\
\bfN: \bfu\in  D(\bfA)\subset H\mapsto \bfN\bfu:=\Big(-2 {\rm P}\,(\bfomega\times\bfv)-{\rm P}\,(\bfv\cdot\nabla\bfv)\,,\ -\bfomega\times\mathbb I\cdot(\bfomega-\bfa)\Big)^\top\in H
\ea
\eeq{1.6}
where $\mathrm P$ is the Helmholtz projection from $L^2(\calc)$ onto $L^2_\sigma(\calc)$. 
As a consequence, the system of equations \eqref{1.3}--\eqref{1.4} can be formally written as the following evolution equation in the Hilbert space $H$
\be
\ode{}t\bfI\bfu+\bfA\bfu+\bfB\bfu=\bfN\bfu\,.
\eeq{1.7}
We shall now state some important properties of the above operators. In the first place, we observe that, by well known results on the Stokes operator, $\bfA_0:=-\nu{\rm P}\Delta$, with
$$
D(\bfA_0)=W^{2,2}(\calc)\cap W_0^{1,2}(\calc)\cap L^2_\sigma(\calc)\,,
$$
it follows that $\bfA$ is  selfadjoint, positive and with compact inverse. For $\alpha\in (0,1)$,  its fractional powers are given by
$$
\bfA^\alpha\bfu:=(\bfA_0^\alpha\bfv\,,\ \bfomega)^\top\,,\ \ \bfA^{-\alpha}:=(\bfA^{-1})^\alpha\,.
$$
Moreover, the operator $\bfB$ is, obviously, linear and bounded. Furthermore, 
in \cite[\S 6.2.3]{kk} it is shown that the bounded operator $\bfI$ is   positive and invertible. In addition, it is easy to check that $\bfI$ is symmetric, namely,
$$
\langle \bfI\bfu_1,\bfu_2\rangle=\langle \bfu_1,\bfI\bfu_2\rangle\,,\ \mbox{$\bfu_1,\bfu_2\in H$}\,,
$$
so that  $\bfI$ is also selfadjoint. Finally, by using classical embedding theorems, one can easily check that $\bfN$ is well defined. As a result, setting
\be
\textsf{\textbf{A}}:=\bfI^{-1}\bfA
\,,\ \ 
\textsf{\textbf{B}}:=\bfI^{-1}\bfB\,,\ \ \textsf{\textbf{N}}:=\bfI^{-1}\bfN\,.
\eeq{1.8}
equation
\eqref{1.7} is equivalent to the following one
\be
\ode\bfu t+\textsf{\textbf{A}}\bfu+\textsf{\textbf{B}}\bfu=\textsf{\textbf{N}}\bfu\ \ \mbox{in $H$}\,.
\eeq{1.9}
From the properties of $\bfI$ and $\bfA$, it readily follows that the operator
 $\textsf{\textbf{A}}$,  with $D(\textsf{\textbf{A}})=D(\bfA)$,  is sectorial \cite[Definition 1.3.1]{Henry} and has a purely discrete, positive  spectrum.\par 
By employing classical methods, it is  now rather straightforward to prove  local in time existence and uniqueness of ``strict" solutions to \eqref{1.9}, corresponding to appropriate initial data. More precisely, we have the following
\Bl Let $\beta\in [\frac43,1)$, and $\bfu_0\in D(\bfA^\beta)$. Then, there exists $t^*\in (0,\infty]$ and a unique function $\bfu=\bfu(t)$ defined for all $t\in [0,t_1]$, $t_1< t^*$, such that
$$
\bfu\in C([0,t_1];D(\bfA^\beta))\cap C((0,t_1];D(\bfA))\cap C^1((0,t_1];H)\,,
$$
solving \eqref{1.9} (or, equivalently, \eqref{1.7}) for all $t\in (0,t_1]$, with $\bfu(0)=\bfu_0$. Moreover, if $t^*<\infty$, then $\|\bfA^\alpha\bfu(t)\|\to \infty$ as $t\to t^*$. 
\EL{1.1}
{\em Proof.} The proof is quite standard (see, e.g., \cite[Lemma 5.1]{KiKi} and \cite[Theorem 3.1 in Chapter 6]{Pazy}) and we shall only sketch it here. Observing that $\textsf{\textbf{A}}$, being sectorial,  
is the generator of the analytic semigroup ${\rm e}^{-{\textsf{\textbf{A}}}\,t}$, 
we begin to consider the following integral equation
\be 
\bfw(t)={\rm e}^{-\textsf{\textbf{A}}\,t}\bfw_0+\int_0^t\textsf{\textbf{A}}^\beta{\rm e}^{-\textsf{\textbf{A}}\,(t-s)}\bfG(\bfw)\,ds\,,\ \ \bfG(\bfw):=-\textsf{\textbf{B}}(\textsf{\textbf{A}}^{-\beta}\bfw)+\textsf{\textbf{N}}(\textsf{\textbf{A}}^{-\beta}\bfw)\,.
\eeq{1.10}
with $\bfw_0=\textsf{\textbf{A}}^{\beta}\bfu_0$. By a classical result  \cite[Lemma 3]{KaFu}
it follows that
\be
\|{\rm P}(\bfv_1\cdot\nabla\bfv_1-\bfv_2\cdot\nabla\bfv_2)\|_2\le c_1\,(\| 
\bfA_0^{\beta}\bfv_1\|_2+\| 
\bfA_0^{\beta}\bfv_2\|_2)\,\|\bfA_0^\beta(\bfv_1-\bfv_2)\|_2.
\eeq{1.10_1}
Also, since $\bfI^{-1}$ is bounded, 
we have (Heinz inequality),
\be
\|\bfA^\alpha\bfu\|\le c_2\,\|\textsf{\textbf{A}}^\alpha\bfu\|\le c_3\,\|\bfA^\alpha\bfu\|\,,\ \ \alpha\in [0,1]\,.
\eeq{EmNa}
From these two properties, we then easily show that
\be
\|\bfG(\bfw_1)-\bfG(\bfw_2)\|\le c_4\|\bfw_1-\bfw_2\|\,,\ \ \mbox{all $\bfw_1,\bfw_2\in H$ in a neighborhood of $\bfw=\bf0$\,.} 
\eeq{lof}
It then follows (e.g., \cite[p. 196--197]{Pazy}) that \eqref{1.10} has one and only one solution $\bfw=\bfw(t)$, $t\in[0,t^*)$ for some $t^*>0$, with $\bfw\in C([0,t_1];H)$, all $t_1\in (0,t^*)$.  As a result, the field $\bfu:=\textsf{\textbf{A}}^{-\beta}\bfw$ is in $C([0,t_1];D(\bfA^\beta))$ and is a mild solution to \eqref{1.9}, namely,
\be
\bfu(t)={\rm e}^{-\textsf{\textbf{A}}\,t}\bfu_0+\int_0^t{\rm e}^{-\textsf{\textbf{A}}\,(t-s)}[-\textsf{\textbf{B}}\bfu+\textsf{\textbf{N}}\bfu]\,ds\,.
\eeq{1.11}
However, again by standard arguments (e.g., \cite[p. 198]{Pazy}),  from \eqref{1.11} one shows that, in fact, $\bfu\in C^1((0,t_1];H)\cap C((0,t_1]; D(\bfA))$, so that $\bfu$ is a ``strict" solution to \eqref{1.9} and is unique. Finally, it is clear that the above procedure can be extended to provide a solution beyond any time $\tau\in [t_1,t^*)$ if $\|\bfA^\beta\bfu(\tau)\|<\infty$, whereas, if $t^*<\infty$, it will fail if and only if $\lim_{t\to t^*}\|\bfA^\beta\bfu(t)\|=\infty$.
\QED

\setcounter{section}{2}\setcounter{equation}{0}
\renewcommand{\theequation}{{2}.\arabic{equation}}
\setcounter{lemm}{0}
\section*{\normalsize 2. The Spectrum of the Linearized Operator}
Let
\be
\bfL:=\textsf{\textbf{A}}+\textsf{\textbf{B}}\,, \ \ D(\bfL)=D(\textsf{\textbf{A}})\equiv D(\bfA).
\eeq{1.12}
The spectral properties of the operator $\bfL$ will play a primary role in establishing  stability of the steady-state motion \eqref{1.0} and the long-time behavior of solutions to the nonlinear problem \eqref{1.9} (or, equivalently, \eqref{1.7}). To this end, we begin to show the following preliminary result.
\Bl The spectrum, $\sigma(\bfL)$, of $\bfL$ 
consists of eigenvalues of finite multiplicity that can cluster only at infinity. Moreover, $0$ is an  eigenvalue with corresponding 
eigenspace 
$$
{\sf N}[\bfL]=\big\{\bfu^{(0)}\in H: \bfu^{(0)}=(\bfv\equiv\0,\bfomega^{(0)})^\top, \ \mbox{for some $\bfomega^{(0)}\in\cals(\lambda)$}\big\}\,,
$$
and 
\be{\dim }\,({\sf N}[\bfL])={\dim }\,(\cals (\lambda))=m\,,\ \ 1\le m\le 3.\eeq{dim}
In addition, the range of $\bfL$, ${\sf R}[\bfL]$, is closed and  the following decomposition holds,
\be
H={\sf N}[\bfL]\oplus {\sf R}[\bfL]\,,
\eeq{nj}
namely, $0$ is a semisimple eigenvalue. 
\EL{1.2}
{\em Proof.} Since $\textsf{\textbf{B}}$ is bounded and $\bfA$ is positive definite, all $\zeta\in\real$ with $\zeta$ sufficiently negative  are in the resolvent set of $\bfL$.  For one of these values of $\zeta$ and for a given $\bfg\in H$, consider the equation
\be
\bfL\bfu-\zeta\,\bfu=\bfg\,.
\eeq{1.13}
Taking into account that ($\bar{\textcolor{white}{\bfu}}$\ :=\ complex conjugate)
$$
\langle \bfA\bfu,\bar{\bfu}\rangle =\nu\|\nabla\bfv\|_2^2+|\bfomega|^2\,,\ \ c_1\|\bfu\|^2\ge \langle\bfI\bfu,\bar{\bfu}\rangle\ge c_0\|\bfu\|^2\,,
$$
by first applying the operator $\bfI$ to both sides  of \eqref{1.13} and then taking the scalar product with $\bar\bfu$, it  follows that
$$
\nu\|\nabla\bfv\|_2^2+|\bfomega|^2-c_0(\zeta+\|\bfB\|)\|\bfu\|^2\le c_1\|\bfg\|\,\|\bfu\|\,.
$$
From the latter expression, by means of the Poincar\'e inequality,  we easily show that, for $\zeta<-\|\bfB\|$, 
$$
\|\bfv\|_{1,2}+|\bfomega|\le c_2\|\bfg\|
$$
which, in turn, by the compact embedding $W^{1,2}(\calc)\subset L^2(\calc)$, implies that $\bfL$ has a compact resolvent. The first part of the lemma is then a consequence of  \cite[Theorem 6.29]{Kato}. To show the second part, we observe that, taking into account \eqref{1.6}, \eqref{1.8} and \eqref{1.12}, the equation
$$
\bfL\bfu=\0\,
$$
is equivalent to the system
$$\ba{ll}\medskip
-\nu\Delta\bfv+2\bfomega_0\times\bfv=-\nabla p\,,\ \ 
\Div\bfv=0\,,\ \mbox{in $\calc$}\,;\ \ \bfv=\0\  \mbox{at $\partial\calc$}\\
\bfomega_0\times\mathbb I\cdot\bfomega+\bfomega\times\mathbb I\cdot\bfomega_0-\bfomega_0\times\mathbb I\cdot\bfa=0
\ea
$$
for some $p\in W^{1,2}(\calc)$. Dot-multiplying both sides of the first equation by ${\bfv}$ and integrating by parts over $\calc$, with the help of the second and third equation we show $\|\nabla\bfv\|_2=0$, which furnishes $\bfv\equiv\0$. Replacing this information back into the last equation, we get
$$
\bfomega_0\times\mathbb I\cdot\bfomega+\bfomega\times\mathbb I\cdot\bfomega_0=\0\,,
$$
which, by \lemmref{1.0}, furnishes $\bfomega\in\cals(\lambda)$.
Therefore, an eigenvector $\bfu^{(0)}$ associated to the  eigenvalue 0 must be of the form $\bfu^{(0)}=(\bfv\equiv\0,\bfomega^{(0)})^\top$, with $\bfomega^{(0)}\in \cals(\lambda)$. This property also implies \eqref{dim}.
We next show that $\bfL$ is an (unbounded) Fredholm operator of index 0. Actually, since $D(\bfA)$ is, obviously, dense in $H$, the same property holds for $D(\bfL)$. Moreover,  by  well-known results concerning the Stokes operator, we have ${\sf N}[\bfA]=\{\0\}$ and ${\sf R}[\bfA]=H$, so that ${\rm ind}(\bfA)=0$. It 
is easy to see that $\bfB$ is $\bfA$-compact, namely, if $\{\bfu_n\}$ is a sequence such that
$$
\|\bfu_n\|+\|\bfA\bfu_n\|\le C
$$
with $C$ independent of $n$, there is $\bfu\in H$ such that
$$
\lim_{n\to\infty}\|\bfB\bfu_n-\bfB\bfu\|=0\,
. 
$$
This property is an elementary consequence of  the classical inequality (e.g. \cite[Theorem IV.6.1]{GaB})
\be
\|\bfA_0\bfv\|_2\le \textcolor{black}{3^{\frac12}}\,\|\bfv\|_{2,2}\le \gamma\, \|\bfA_0\bfv\|_2\,,\ \ \gamma=\gamma(\calc)>0\,,
\eeq{EmNa_0}
\textcolor{black}{the compact embedding} $W^{2,2}\subset L^2$,  and \textcolor{black}{the definition of} $\bfB$ given in \eqref{1.6}$_3$. Therefore, $\bfA+\bfB$ is Fredholm of index 0. Since $\bfI^{-1}$ is \textcolor{black}{a} homeomorphism, by the product property, we find that also $\bfL$ is Fredholm of index 0.
Now, assume ${\rm dim}\,\cals(\lambda)=1$. By what we have already shown, ${\rm dim}({\sf N}[\bfL])=1$ as well, which, by the Fredholm property of $\bfL$, implies ${\rm codim}({\sf R}[\bfL])=1$. Let us prove  that 
\be
{\sf R}[\bfL]\cap {\sf N}[\bfL]=\{\0\}\,. 
\eeq{df}
In fact, suppose there is $\bfu^{(0)}=(\0,\bfomega^{(0)})^\top\in {\sf N}[\bfL]$ such that
$$
\bfL\bfu=\bfu^{(0)}\,,
$$
or, equivalently,
$$
(\bfA+\bfB)\bfu=\bfI\bfu^{(0)}\,.
$$
Recalling \eqref{1.6},
 this equation can be rewritten as follows
\be
\ba{ll}\medskip
-\nu\Delta\bfv+2\bfomega_0\times\bfv=-\nabla p+\bfomega^{(0)}\times\bfx\,,\ \ 
\Div\bfv=0\,,\ \mbox{in $\calc$}\,;\ \ \bfv=\0\  \mbox{at $\partial\calc$}\\
\bfomega_0\times\mathbb I\cdot\bfomega+\bfomega\times\mathbb I\cdot\bfomega_0-\bfomega_0\times\mathbb I\cdot\bfa=\mathbb I\cdot\bfomega^{(0)}\,.
\ea
\eeq{cu1}
Taking into account that $\mathbb I\cdot\bfomega_0=\lambda\,\bfomega_0$ and that, by assumption, $\bfomega^{(0)}=\alpha\,\bfomega_0$, $\alpha\in\real$, from \eqref{cu1} we get, in particular,
$$
\bfomega_0\times\big(\mathbb I\cdot\bfomega-\lambda\,\bfomega-\mathbb I\cdot\bfa\big)=\lambda\,\alpha\,\bfomega_{0}\,,
$$
which implies $\bfomega^{(0)}=\0$, namely, $\bfu^{(0)}=\0$, and as a result we conclude \eqref{df}. Likewise, if ${\rm dim}\,(\cals(\lambda))\equiv{\rm dim}\,({\sf N}[\bfL])=3$, we have $(\mathbb I)_{ij}=\lambda\, \delta_{ij}$, so that  the last equation in \eqref{cu1} furnishes
\be
\bfomega^{(0)}\cdot\bfa=0\,.
\eeq{si}
On the other side, by dot-multiplying both sides of the first equation in \eqref{cu1} by $\bfv$,  integrating by parts over $\calc$, and recalling \eqref{1.5}, we show
$$
\nu\|\nabla\bfv\|_2^2=\bfomega^{(0)}\cdot\bfa
$$ 
which in view of \eqref{si} implies $\bfv\equiv\0$. Replacing this information into \eqref{cu1}$_1$ entails
\be
\bfomega^{(0)}\times\bfx=\nabla p\,,
\eeq{cp}
so that, by operating with $\curl$ on both sides, we conclude, also in this case, $\bfu^{(0)}=0$ and the validity of \eqref{df}.
Let us, finally, consider the case $m=2$. 
In this regard, we begin to notice that from\eqref{cu1}$_4$
 it follows  $\bfomega_0\cdot\bfomega^{(0)}=0$.  To  
 fix the ideas, we assume $\lambda\equiv A=B$ $\textcolor{black}{(<C)}$, the case $B=C$ $(>A)$ being treated in an entirely  similar way. Thus, taking (for instance)
\be
\bfomega_0=\omega_0\,\bfe_1\,, \ \bfomega^{(0)}=\omega^{(0)}\,\bfe_2\,,\eeq{Nash0} equation   \eqref{cu1}$_4$ becomes
$$
\bfomega_0\times [(C-A)\bfomega]=\bfomega_0\times\mathbb I\cdot\bfa+A\,\bfomega^{(0)}:=\bfF\,,
$$
where $\bfomega=\omega\,\bfe_3$. Since $\bfF\cdot\bfomega_0=0$, this equation is solvable for $\bfomega$ and we get 
\be
\bfomega=\frac{\bfF\times\bfomega_0}{\omega_0^2\,(C-A)}\,.
\eeq{Nas}
Observing that
$$
\bfF\times\bfomega_0=\omega_0^2\,\mathbb I\cdot\bfa-(\bfomega_0\cdot\mathbb I\cdot\bfa)\,\bfomega_0+A\,\bfomega^{(0)}\times\bfomega_0\,,
$$
from \eqref{Nash0} and \eqref{Nas} we deduce that the  condition 
$\bfomega\cdot\bfomega^{(0)}=0$ implies
\be
\bfomega^{(0)}\cdot\mathbb I\cdot\bfa=0\,.
\eeq{Nash1}
However, if we dot-multiply both sides of \eqref{cu1}$_1$ by $\bfv$, integrate by parts over $\calc$ and use \eqref{cu1}$_{2,3}$ and \eqref{1.5}, we get
$$
\|\nabla\bfv\|_2^2=\bfomega^{(0)}\cdot\mathbb I\cdot\bfa, 
$$
which, combined with \eqref{Nash1} furnishes $\bfv\equiv\0$. Replacing the latter in \eqref{cu1}$_1$ entails again \eqref{cp} that, as shown previously, implies $\bfomega^{(0)}=\0$ and  
 the proof of \eqref{df} is completed.  
Now, since ${\rm codim}({\sf R}[\bfL])=m$, there exists at least one $S\textcolor{black}{\subset} H$ such that $H=S\oplus{\sf R}[\bfL]$, with $S\cap {\sf R}[\bfL]=\{\0\}$. However, ${\rm dim}\,(S)={\rm dim}\,({\sf N}[\bfL])=m$ and \eqref{df} holds, so that we may take $S={\sf N}[\bfL]$, as claimed. 

\QED 
\par
Let $\calq$ and $\calp$ be the spectral projections according to the spectral sets
$$\sigma_0(\bfL):=\{0\}\,,\ \
\sigma_1(\bfL):= \sigma(\bfL)\backslash\sigma_0(\bfL)\,,
$$
which are well defined in view of \lemmref{1.2}.
Thus, setting
\be
H_0:=\calq(H)\,,\ \ H_1:=\calp(H)\,,
\eeq{Sab}
we have the following decomposition
\be
H=H_0\oplus H_1
\eeq{dec}
that completely reduces $\bfL$ into $\bfL=\bfL_0\oplus\bfL_1$
with
\be
\bfL_0:=\calq\bfL=\bfL\calq\,,\ \ \bfL_1:=\calp\bfL=\bfL\calp\,, 
\eeq{1.18_0}
and $\sigma(\bfL_0)\equiv\sigma_0(\bfL)$, $\sigma(\bfL_1)\equiv\sigma_1(\bfL)$ (e.g., \cite[Theorems 5.7-A,B]{Tay}). 
\medskip\par
A fundamental issue in the proof of the nonlinear results that we shall present later on, is the identification of the subspace $H_0$ with ${\sf N}[\bfL]$. In general, this is not true and we only have  ${\sf N}[\bfL]\subseteq H_0$, whereas ${\sf R}[\bfL]\supseteq H_1$. However, if (and only if) the decomposition \eqref{nj} holds, then this property is valid and the above subspaces coincide (e.g. \cite[Proposition A.2.2]{Lun}). Thus, \lemmref{1.2} implies the next one.
\Bl The following characterization holds
$$
H_0={\sf N}[\bfL]\,,\ \ H_1={\sf R}[\bfL]\,.
$$
\EL{2.2}
\Br It is worth noticing that one can show that the ortho-complement of ${\sf N}[\bfL]$ does {\em not} coincide with ${\sf R}[\bfL]$. Therefore, under the assumption of \lemmref{2.2}, the decomposition \eqref{dec} is not orthogonal.  
\ER{Li}

The next result provides a complete  characterization of the distribution of the eigenvalues of $\bfL_1$ in terms of the central moments of inertia of $\mathscr S$ and of the axis where the permanent rotation occurs. 

\Bp Let ${\sf s}_0$ be given by \eqref{1.0}. Then, 
\be
\sigma_1(\bfL)\cap \{{\rm i}\,\mathbb R\}=\emptyset.
\eeq{Lun}
Moreover, we have
\be
\Re[\sigma_1(\bfL)]\subset (0,\infty)\,.
\eeq{1.15}
whenever at least one  of the following conditions holds. 
\begin{itemize}
  \item [{\rm (i)}] $A=B=C$, arbitrary $\bfe$;
  \item [{\rm (ii)}] $A\le B<C$, $\bfe\equiv\bfe_3$\,;
    \item [{\rm (iii)}] $A<B=C$, $\bfe=\gamma_1\bfe_2+\gamma_2\bfe_3$, with $\bfgamma\equiv(\gamma_1,\gamma_2)\in S^1$\,. 
\end{itemize}
\par
Conversely, if  any of the following conditions are met
\begin{itemize}
\item[{\rm (iv)}] $A<B\le C$, $\bfe\equiv\bfe_1$\,;
\item[{\rm (v)}] $A<B<C$, $\bfe\equiv\bfe_2$\,;
\item[{\rm (vi)}] $A=B<C$, $\bfe\equiv\gamma_1\bfe_1+\bfgamma_2\bfe_2$ with $\bfgamma=(\gamma_1,\gamma_2)\in S^1$\,.
\end{itemize}
Then \be
\Re[\sigma_1(\bfL)]\cap (-\infty,0)\neq\emptyset\,.
\eeq{lava}
\par 
\EP{1.2}
{\em Proof.} Consider the linear  problem
\be
\ode{\bfu}t+\bfL\bfu=\0\,,\ \ \bfu(0)=\bfu_0\in H\,.
\eeq{1.37}
We begin to notice that, being a bounded perturbation of the operator $\textsf{\textbf{A}}$, the operator $\bfL$  is the generator of an analytic semigroup (e.g., \cite[Theorem 2.1 in Chapter 3]{Pazy}). As a consequence, the solution to \eqref{1.37} is unique and smooth for all $t\in (0,\infty)$. In particular, the map 
$$
\bfu_0\mapsto \bfu(t;\bfu_0)
$$
with $\bfu(\,\cdot \,;\bfu_0)$ solution to \eqref{1.37} corresponding to the initial data $\bfu_0$, defines a dynamical system in $H$.
Now, the abstract equation \eqref{1.37} is equivalent to the following system
\be
\ba{cc}\medskip
\left.\ba{ll}\medskip
\bfv_t+\dot{\bfomega}\times\bfx+2\bfomega_0\times\bfv-\nu\Delta\bfv-\nabla p=\0\\
\Div\bfv=0\ea\right\}\ \ \mbox{in $\calc\times(0,\infty)$}\,,\\ \medskip
\bfv(x,t)=\0\ \ \mbox{at $\partial\Omega$}\,,\\
\mathbb I\cdot(\dot{\bfomega}-\dot{\bfa})+\bfomega_0\times\mathbb I\cdot\bfomega+\bfomega\times\mathbb I\cdot\bfomega_0-\bfomega_0\times\mathbb I\cdot\bfa=\0\,,
\ea
\eeq{1.41}
where, we recall,  $\bfa$ is defined in \eqref{1.5}. 
Next, 
for a given vector ${\bfh}\in \real ^3$, we  write ${\bfh}={\bfh}_{\perp}+\bfh_{\|}$, where  $\bfh_{\|}= h_{\|}\,\bfe$. Clearly, $\bfh_\perp\in \cals(\lambda)^\perp$ if ${\rm dim}\,\cals(\lambda)=1$, while $\bfh_\perp=\0$ if ${\rm dim}\,\cals(\lambda)=3$.
Since
$$
\bfomega_0\times\mathbb I\cdot\bfomega_{\|}+\bfomega_{\|}\times\mathbb I\cdot\bfomega_0=\0,
$$
from \eqref{1.41}$_4$, we  deduce
\be\ba{ll}\medskip
\mathbb I\cdot(\dot{\bfomega}_{\perp}-\dot{\bfa}_{\perp})+\bfomega_0\times\mathbb I\cdot\bfomega_{\perp}+\bfomega_{\perp}\times\mathbb I\cdot\bfomega_0-\bfomega_0\times\mathbb I\cdot\bfa_{\perp}=\0
\\ 
\dot{\omega}_{\|}=\dot{a}_{\|}\,.
\ea
\eeq{1.42}
Setting $\boms:=\bfomega_{\perp}-\bfa_{\perp}$ and taking into account \eqref{1.42}$_2$ we can then rewrite \eqref{1.41} in the following equivalent way:
\be
\ba{cll}\medskip
\left.\ba{ll}\medskip
\bfv_t+\dot{\bfa}\times\bfx+\dot{\bfomega}_*\times\bfx+2\bfomega_0\times\bfv-\nu\Delta\bfv-\nabla p=\0\\
\Div\bfv=0\ea\right\}\ \ \mbox{in $\calc\times(0,\infty)$}\,,\\ \medskip
\bfv(x,t)=\0\ \ \mbox{at $\partial\Omega$}\,,\\ \medskip
\mathbb I\cdot\dot{\bfomega}_*+\bfomega_0\times\mathbb I\cdot\boms+(\boms+\bfa_{\perp})\times\mathbb I\cdot\bfomega_0=\0\,,
\\
\dot{\omega}_{\|}=\dot{a}_{\|}\,.
\ea
\eeq{1.43}
If we dot-multiply \eqref{1.43}$_1$ by $\bfv$, integrate by parts over $\calc$ and take into account \eqref{1.43}$_{2,3}$, we obtain
\be
\ode{E}t+\nu\|\nabla\bfv\|_2^2=\dot{\bfomega}_*\cdot\mathbb I\cdot\bfa_{\perp}\,, 
\eeq{1.44}
where
\be
E:=\frac12 (\|\bfv\|^2_2-\bfa\cdot\mathbb I\cdot\bfa)\,.
\eeq{ene}
Due to \cite[\S\S\, 7.23, 7.2.4]{kk1}, we know that there is $c_0\in (0,1)$ such that
\be
c_0\|\bfv\|_2^2\le 2E\le \|\bfv\|_2^2\,.
\eeq{1.45}
We next dot-multiply \eqref{1.43}$_4$ one time by $\boms$ and a second time by $\bfa_{\perp}$ to get
$$\ba{ll}\medskip
\boms\cdot\mathbb I\cdot\dot{\bfomega}_*+\bfomega_0\times\mathbb I\cdot\boms\cdot\boms+(\boms+\bfa_{\perp})\times\mathbb I\cdot\bfomega_0\cdot\boms=\0
\\
\dot{\bfomega}_*\cdot\mathbb I\cdot\bfa_{\perp}+\bfomega_0\times\mathbb I\cdot\boms\cdot\bfa_{\perp}+(\boms+\bfa_{\perp})\times\mathbb I\cdot\bfomega_0\cdot\bfa_{\perp}=\0\,,
\ea
$$
from which we show
$$
\dot{\bfomega}_*\cdot\mathbb I\cdot\bfa_{\perp}=-\boms\cdot\mathbb I\cdot\dot{\bfomega}_*-\bfomega_0\times\mathbb I\cdot\boms\cdot(\boms+\bfa_{\perp})\,.
$$
Replacing the latter into \eqref{1.44} produces
\be
\ode{}t\big(E+\half \boms\cdot\mathbb I\cdot\boms\big)+\nu\|\nabla\bfv\|_2^2=-\bfomega_0\times\mathbb I\cdot\boms\cdot(\boms+\bfa_{\perp})\,.
\eeq{1.46}
On the other hand, if we dot-multiply both sides of \eqref{1.43}$_4$ by $\mathbb I\cdot\boms$ and recall that $\mathbb I\cdot\bfomega_0=\lambda\,\bfomega_0$, for some $\lambda\in \{A,B,C\}$, we show
$$
\frac12\ode{}t(\mathbb I\cdot\boms)^2=-\lambda \,\bfomega_0\times\mathbb I\cdot(\boms+\bfa_{\perp})\,.
$$
The latter, in conjunction with \eqref{1.46} allows us to conclude
\be
\ode{}t\big(2E+\boms\cdot\mathbb I\cdot\boms-\frac1\lambda\,(\mathbb I\cdot\boms)^2 \big)+2\nu\|\nabla\bfv\|_2^2=0\,.
\eeq{1.47}
This equation is the fundamental tool in our proof, In fact, let us begin to show \eqref{Lun}. Assuming the contrary would imply that \eqref{1.37} has (at least) one non-trivial, smooth  solution $\bfu=(\bfv,\bfomega)^\top$  such that
\be
\bfu(T)=\bfu(0)\,,\ \ \int_0^T\bfu(t)=0\,,\ \ \mbox{some $T>0$}
. 
\eeq{kc}
Integrating both sides of \eqref{1.47} over $[0,T]$ and using \eqref{kc} would then produce
$$
\int_0^T\|\nabla\bfv(t)\|_2^2=0\,,
$$
namely, \be\bfv\equiv\0 \ \ \mbox{in $[0,T]$.}\eeq{kc2} 
If we replace this information back in \eqref{1.43}, we get, in particular,
\be
\dot{\bfomega}_*\times\bfx=\nabla p\,,\ \ 
\mathbb I\cdot\dot{\bfomega}_*+\bfomega_0\times\mathbb I\cdot\boms+\boms\times\mathbb I\cdot\bfomega_0=\0\,.
\eeq{kc1}
Operating with $\curl$ on both sides of the first of these equations we deduce $\boms=\textrm{const.}$. 
This condition combined with \eqref{kc}, \eqref{kc2} and  \eqref{1.43}$_5$  implies 
$\bfomega(t)\equiv\0$, so that the latter and, again, \eqref{kc2} provide that there is no non-trivial solution $\bfu$ to \eqref{1.37} satisfying \eqref{kc}, and \eqref{Lun} is thus established.  With the help of this result and \lemmref{1.2} we then infer
that $\sigma_1(\bfL)$ must satisfy either \eqref{1.15} or \eqref{lava}.
 We next prove the second property stated in the proposition. To this end, we observe that
if \eqref{1.15} is not true, then there exists at least one solution to \eqref{1.37} that becomes unbounded in $H$ as $t\to\infty$.
Thus, 
in order to show our stability claim
it
is sufficient to show that, under any of the assumptions (i)--(iii), all solutions to \eqref{1.37} are bounded. In turn, in view of \eqref{1.43}$_5$, this amounts to show that there is a constant $M>0$, depending on the data, such that
\be
\|\bfv(t)\|_2+|\bfomega_\perp(t)|\le\, M \,,\ \ \mbox{all $t\ge0$.}
\eeq{bu}
Set
$$
G:=\boms\cdot\mathbb I\cdot\boms-\frac1\lambda\,(\mathbb I\cdot\boms)^2
$$
and consider the three different cases (i)--(iii) stated in the proposition.
If $A=B=C$, we get 
$G=0$, so that from \eqref{1.45}, \eqref{1.47} and Poincar\'e inequality, we deduce
$$
\ode{E}t\le -\sigma E\,,\ \sigma=\textrm{const.}>0
$$
that in turn, again by \eqref{1.45}, implies that $\|\bfv(t)\|_2$ is uniformly bounded in time. Since in this case $\bfomega_{\perp}=\0$, \eqref{bu} is proved, and, with it, the proposition in case (i). In case (ii), we get 
$$
G=A\omega_{*1}^2+B\omega_{*2}^2-\frac1C\big(A^2\omega_{*1}^2+B^2\omega_{*2}^2\big)=\alpha\, \omega_{*1}^2+\beta\,\omega_{*2}^2\,,
$$
where
$$
 \alpha:=\frac1C A(C-A)\,,\  \beta:=\frac1C B(C-B) \, >0\,.
$$
So, under the given assumptions, $G$ is a positive definite quadratic form in the components of $\bfomega_*$. This property along with \eqref{1.45} and \eqref{1.47} implies that there exists a constant $M>0$ such that
$$
\|\bfv(t)\|_2^2+|\bfomega_*(t)|^2\le M\,,\ \ \mbox{ all $t\ge0$,}
$$
which proves condition \eqref{bu}. We next consider the case (iii). Without loss of generality (we can always rotate $\bfe_2$ and $\bfe_3$ appropriately) we take $\bfe\equiv\bfe_3$. We thus deduce
$$
G=A\omega_{*1}^2+B\omega_{*2}^2-\frac1B(A^2\omega_{*1}^2+B^2\omega_{*2}^2)=\frac AB(B-A)\omega_{*1}^2\,,
$$
which with the help of  \eqref{1.45} and \eqref{1.47} entails, for some constant $M>0$,
\be
\|\bfv(t)\|_2+|\omega_1(t)|\le M\,,\ \ \mbox{all $t\ge 0$}\,.
\eeq{gp}
However, taking the projection of  \eqref{1.37} along the subspace $H_0$ and recalling \lemmref{2.2}, we also obtain $\omega_2(t)=\textrm{const.}$, so that \eqref{bu} follows from the latter and \eqref{gp}. 
 Finally, we prove the last claim in the proposition. If \eqref{lava} were not true, then all solutions to \eqref{1.37} must be uniformly bounded in time, so that, in particular, they must satisfy \eqref{bu}. We shall then show that, in such a case, the following relation holds
\be 
2\nu\int_0^\infty\|\nabla\bfv(s)\|_2^2ds=2E(0)+G(0)\,.
\eeq{1.54}
To this end, by integrating both sides of  \eqref{1.47} and using \eqref{bu} we deduce that  
$$
\int_0^\infty\|\nabla\bfv(t)\|_2^2<\infty\,,
$$
which, by Poincar\'e inequality, entails 
\be
\int_0^\infty\|\bfv(t)\|_2^2dt\le M\,.
\eeq{1.49}
Now, \eqref{bu} 
along with \eqref{1.43}$_4$, furnishes  
$$
|\dot{\bfomega}_*(t)|\le M_1\ \ \mbox{ all $t>0$,}
$$
for another constant $M_1>0$. Replacing this information on the right-hand side of \eqref{1.44} and using \eqref{1.45} and Schwarz and Poincar\'e inequalities, we show
$$
\ode{E}t+c_1\, E\le c_2\,E^{\frac12}\,,
$$
which, by a generalized form of Gronwall's lemma \cite[Lemma 2.1]{GMZ} and \eqref{1.45}, \eqref{1.49} implies
\be
\lim_{t\to\infty}\|\bfv(t)\|_2=0\,.
\eeq{1.50}
From \eqref{1.43}$_5$, \eqref{bu} and \eqref{1.50} we infer that the orbits generated by the solutions to \eqref{1.37} through any initial data $\bfu_0$ are compact and, therefore, the $\omega$-limit set is not empty and, in particular, invariant. By \eqref{1.50}, $\bfv\equiv \0$ on this set so that by taking the $\curl$ of both sides of  \eqref{1.43}$_1$ (with $\bfv\equiv\0$) we derive $\boms=\bar{\bfomega}=\textrm{const.}$, which once replaced in \eqref{1.43}$_4$ entails that $\bar{\bfomega}$ must satisfy
$$
\bfomega_0\times\mathbb I\cdot\bar{\bfomega}+\bar{\bfomega}\times\mathbb I\cdot\bfomega_0=\0\,.
$$
From  \lemmref{1.0} and the latter we thus deduce 
\be
\bar{\bfomega}\in \cals(\lambda).\eeq{jb} At this point, let us first discuss  the cases (iv) and (v) stated in the proposition. Then, in both situations,  $\bar{\bfomega}$ must also belong to $\cals(\lambda)^\perp$, which produces $\bar{\bfomega}=\0$. Therefore, by definition of $\omega$-limit set, we conclude
\be
\lim_{t\to\infty}|\boms(t)|=0\,.
\eeq{ms}
If we now integrate both sides of \eqref{1.47} from $0$ to $t$ and then let $t\to\infty$, with the help of \eqref{1.50}, \eqref{1.45} and \eqref{ms} we show \eqref{1.54}. 
Now, 
assume condition (iv). In that case, we have $\boms=\omega_{*2}\bfe_2+\omega_{*3}\bfe_3$. Therefore
\be
G=\frac{B}A(A-B)\,\omega_{*2}^2+\frac{C}A(A-C)\,\omega_{*3}^2
\eeq{1.53}
Thus, from \eqref{1.53} and \eqref{1.54}, we deduce
$$
0\le 2E(0)+\frac{B}A(A-B)\,\omega_{*2}^2(0)+\frac{C}A(A-C)\,\omega_{*3}^2(0)
$$
which, in view of the assumptions on $A,B,$ and $C$ cannot be true if we pick initial data such that
$$
2E(0)<\frac{B}A(B-A)\,\omega_{*2}^2(0)+\frac{C}A(C-A)\,\omega_{*3}^2(0)\,. 
$$
Thus, \eqref{lava} is established. Likewise, in the case (v), we have
$$
G=\frac{A}B(B-A)\,\omega_{*1}^2+\frac{C}B(B-C)\,\omega_{*3}^2\,,
$$
and, arguing  as before, from \eqref{1.54} we deduce
$$
0\le 2E(0)+\frac{A}B(B-A)\,\omega_{*1}^2(0)+\frac{C}B(B-C)\,\omega_{*3}^2(0)\,,
$$ 
which cannot hold if we choose (for example) 
$$
2E(0)<\frac{C}B(C-B)\,\omega_{*3}^2(0)\,,\ \ \omega_{*1}(0)=0\,.
$$
We now turn to the case (vi), and begin to show that also in this case \eqref{1.54} holds. To this end, without loss of generality (it is enough to rotate $\bfe_1$ and $\bfe_2$ appropriately) we assume $\bfe=\bfe_1$. Therefore, we have 
\be
G= A\omega_{*2}^2+C\omega_{*3}^2-\frac1A(A^2\omega_{*2}^2+C^2\omega_{*3}^2)=\frac{C}A(A-C)\omega_{*3}^2\,.
\eeq{jb1}
On the other hand, from \eqref{jb} we deduce 
$\bar\omega_3=0$, which gives
\be
\lim_{t\to\infty}\omega_{*3}(t)=0\,.
\eeq{jb2}
Thus, integrating both sides of \eqref{1.47} from $0$ to $t$,  letting $t\to\infty$, and using \eqref{jb1} and \eqref{jb2} we establish \eqref{1.54} also in the case (vi).
As a consequence,  from \eqref{1.54} we obtain a contradiction if we choose initial data such that
$$
E(0)<\frac{C}A(C-A)\omega_{*3}^2(0)\,.
$$
The proof of the proposition is completed. 
\QED
\Br In physical terms, \propref{1.2}  may be restated by saying that a permanent rotation is ``linearly stable"  if and only if it occurs around an axis of maximum moment of inertia. This is in total agreement with Kelvin's experiment.  
\ER{nj2}
\setcounter{section}{3}\setcounter{equation}{0}
\renewcommand{\theequation}{{3}.\arabic{equation}}
\setcounter{lemm}{0}
\section*{\normalsize 3. Global Existence and  Stability of Permanent Rotations. A Full Explanation of Kelvin's Experiment
}
The main objective of this section is to use the results established in the previous one, to prove necessary and sufficient conditions for the stability of the permanent rotation \eqref{1.0}.
To this purpose, we give the following definition.
\Bd The permanent rotation ${\sf s}_0$ in \eqref{1.0} is called {\em stable} if for any $\varepsilon>0$ there is $\delta=\delta(\varepsilon)>0$ such that 
$$
\|\bfA^\beta\bfu_0\|<\delta\ \ \Longrightarrow\ \ \sup_{t\ge 0}\|\bfA^\beta\bfu(t)\|<\varepsilon
$$
for some $\beta\in [0,1]$ and all solutions $\bfu=\bfu(t)$ to \eqref{1.7} with $\bfu(0)=\bfu_0$. Also,  ${\sf s}_0$ is called {\em unstable} if it is not stable. Furthermore, ${\sf s}_0$ is  {\em asymptotically stable}, if it is stable and there exist $\rho_0>0$ such that
$$
\|\bfA^{\beta}\bfu_0\|<\rho_0\\ \ \Longrightarrow\ \ \lim_{t\to\infty}\|\bfA^{\beta}\bfu(t)-\bfu^{(0)}\|=0\,,
$$
for some  $\bfu^{(0)}=(\0,\bfomega^{(0)})^\top\in H$ with $\bfomega^{(0)}\in\cals(\lambda)$.
Finally, ${\sf s}_0$ is {\em exponentially stable} if it is asymptotically stable and there are constants $C,\kappa>0$ such that
$$
\|\bfA^\beta\bfu(t)-\bfu^{(0)}\|\le C\,\|\bfA^\beta\bfu(0)\|\,{\rm e}^{-\kappa\,t}\,,\ \ \mbox{all $t\ge 0$.}
$$
\EED{1.1}
\Br Because of conservation of total angular momentum for $\mathscr S$, in the asymptotic stability definition, we cannot expect $\bfu^{(0)}=\0$, due to the fact that ${\sf s}_0$ in \eqref{1.0} is non-trivial. As a matter of fact, $\bfu^{(0)}=\0$ if and only if the initial data $\bfu_0:=(\bfv(0),\bfomega(0))^\top$ satisfy  \cite[Remark 5]{DGMZ}
$$
\mathbb I\cdot(\bfomega(0)+\bfomega_0)+\int_\calc\bfx\times\bfv(0)=\0\,,
$$
a non-generic condition that is physically irrelevant. 
\ER{3.1}
\par
The next result provides a suitable ``linearization principle" for the equation \eqref{1.7}.
\Bp \textcolor{black}{Let ${\sf s}_0$ be given by \eqref{1.0}.} The following stability properties are valid.
\begin{itemize}
\item[{\rm (a)}] 
If 
\eqref{1.15} holds,
then  there exists $\gamma_0>0$ such that if for some $\beta\in [\frac34,1)$
$$
\|\bfA^{\beta}\bfu_0\|<\gamma_0\,,
$$
the unique solution $\bfu$ to \eqref{1.9} (or, equivalently, \eqref{1.7}) constructed in \lemmref{1.1} exists for all $t>0$ (namely, we can take $t^*=\infty$). 
Moreover, ${\sf s}_0$ is exponentially stable.  
\item[{\rm (b)}] Conversely, if  \eqref{lava} holds 
then ${\sf s}_0$ is unstable.
\end{itemize}
\EP{1.1}
{\em Proof.} In view of  \lemmref{2.2}, for any $\bfu:=(\bfv,\bfomega)^\top\in H$, the spectral projections $\calq,\calp$ satisfy
$$
\calq\bfu:=\bfu^{(0)}\equiv(\0,\bfomega^{(0)})^\top\in {\sf N}\,[\bfL]\equiv H_0\,,\ \ \calp\bfu:=\bfu^{(1)}\equiv(\bfv,\bfomega^{(1)})^\top\in {\sf R}\,[\bfL]\equiv H_1\,. 
$$
Applying  $\calq$ and $\calp$  on both sides of \eqref{1.9} 
 and taking into account \eqref{1.18_0} we easily show 
$$\ba{rl}\medskip
\ode{\bfu^{(1)}}t+\bfL_1\bfu^{(1)}=&\!\!\!\mathcal P\bfI^{-1}\bfM(\bfu^{(1)},\bfu^{(0)})\\
\ode{\bfu^{(0)}}t +\bfL_0\bfu^{(0)}=&\!\!\!\mathcal Q\bfI^{-1}\bfM(\bfu^{(1)},\bfu^{(0)})\,,
\ea
$$
where  
\be \bfM(\bfu^{(1)},\bfu^{(0)}):=\bfN(\bfu^{(1)}+\bfu^{(0)})\,.\eeq{NaRe} 
Moreover, observing that $\bfL_0\bfu^{(0)}=\0$, the previous equations simplify to the following ones
\be\ba{rl}\medskip
\ode{\bfu^{(1)}}t+\bfL_1\bfu^{(1)}=&\!\!\!\mathcal P\bfI^{-1}\bfM(\bfu^{(1)},\bfu^{(0)})\\
\ode{\bfu^{(0)}}t =&\!\!\!\mathcal Q\bfI^{-1}\bfM(\bfu^{(1)},\bfu^{(0)})\,.
\ea
\eeq{1.20}
In view of \lemmref{1.2} we know that \eqref{1.20} has one and only one  solution --in the class specified there-- in the time interval $[0,t^*)$, which  can be extended to a global solution provided we show the existence of $\rho>0$ such that
\be 
\sup_{t\in[0,t^*)}\|\bfA^\beta\bfu(t)\|\le \rho\,.
\eeq{1.21}
In this regard, we begin to observe that  since the operator $\bfL$ is the generator of an analytic semigroup in $H$,  so is $\bfL_1$ in $H_1$.  Thus, for all $t\in[0,t^*)$ from \eqref{1.20}$_1$ we have
\be
\bfu^{(1)}(t)={\rm e}^{-\bfL_1 t}\bfu_0^{(1)}+\int_0^t{\rm e}^{-\bfL_1 (t-s)}[\mathcal P\bfI^{-1}\bfM(\bfu^{(1)}(s),\bfu^{(0)}(s))]ds\,.
\eeq{1.22}
Also, by assumption and \lemmref{1.2},  there is $\gamma>0$ such that 
\be
\Re[\sigma(\bfL_1)]>\gamma>0\,,  
\eeq{1.17}
which implies that the fractional powers $\bfL_1^\alpha$, $\alpha\in (0,1)$, are well defined in $H_1$. Thus, setting
$$
\bfw:={\rm e}^{bt}\bfL_1^\beta\bfu^{(1)}\,,\ \ 0<b<\gamma\,,
$$
from \eqref{1.22} we get
\be
\bfw(t)={\rm e}^{bt}{\rm e}^{-\bfL_1 t}\bfL_1^\beta\bfu_0^{(1)}+\int_0^t{\rm e}^{bt}\bfL_1^{\beta} {\rm e}^{-\bfL_1 (t-s)}[\mathcal P\bfI^{-1}\bfM({\rm e}^{-bs}\bfL_1^{-\beta}\bfw(s),\bfu^{(0)}(s))]ds\,.
\eeq{1.23}
We now 
make the obvious but {\em crucial} observation that 
$$
\bfomega^{(0)}\times\mathbb I\cdot\bfomega^{(0)}=\0\,. 
$$
So,  from \eqref{NaRe} and \eqref{1.6}$_4$, we obtain
\be\ba{rl}\medskip
\bfM(\bfu^{(1)},\bfu^{(0)})=
\Big(-&\!\!\!\!\!\! {\rm P}\,[2(\bfomega^{(1)}+\bfomega^{(0)})\times\bfv+\bfv\cdot\nabla\bfv]\,,\\
&\! -\bfomega^{(1)}\times\mathbb I\cdot(\bfomega^{(0)}+\bfomega^{(1)})-\bfomega^{(0)}\times\mathbb I\cdot\bfomega^{(1)}+(\bfomega^{(0)}+\bfomega^{(1)})\times\mathbb I\cdot\bfa
\Big)^\top\,.\ea
\eeq{1.26}
With the help of \eqref{1.10_1}, we show
\be\ba{rl}\medskip
\|\bfM(\bfu^{(1)},\bfu^{(0)})\|&\!\!\!\!\le c\,\big[(|\bfomega^{(1)}|+|\bfomega^{(0)}|)(\|\bfv\|_2 +|\bfomega^{(1)}|)+\|\bfA_0^\beta\bfv\|^2\big]\\
&\le c\,\big[(\|\bfu^{(1)}\|+\|\bfu^{(0)}\|)\|\bfu^{(1)}\|+\|\hat{\bfA}^\beta\bfu^{(1)}\|^2\big]\,,
\ea
\eeq{1.27}
where $\hat{\bfT}$ denotes the restriction of the operator $\bfT$ to $H_1$. Since $\textsf{\textbf{B}}$ is a bounded operator, from the definition \eqref{1.12} we have (e.g. \cite[Theorem 1.4.6]{Henry}) 
$$
\|\hat{\textsf{\textbf A}}^\alpha\bfu\|\le c_1\,\|\bfL_1^\alpha\bfu\|\le c_2\,\|\hat{\textsf{\textbf A}}^\alpha\bfu\|\,,\ \ \alpha\in [0,1]\,,
$$
which, by \eqref{EmNa}, implies
\be
\|\hat{\bfA}^\alpha\bfu\|\le c_1\,\|\bfL_1^\alpha\bfu\|\le c_2\,\|\hat{\bfA}^\alpha\bfu\|\,,\ \ \alpha\in [0,1]\,.
\eeq{1.28}
Consequently, from \eqref{1.27} and the latter we derive
\be
\|\mathcal P\bfI^{-1}\bfM({\rm e}^{-bt}\bfL_1^{-\beta}\bfw(s),\bfu^{(0)}(s))\|\le c_3
\,\big[(\|\bfw\|+\|\bfu^{(0)}\|)\|\bfw\|+\|\bfw\|^2\big]\,.
\eeq{1.29}
Next, observing that, in $H_1$, it is $\|\bfL_1^{\alpha}{\rm e}^{-\bfL_1t}\|\le t^{-\alpha}{\rm e}^{-\gamma t}$, $\alpha\in [0,1]$, from  \eqref{1.23}, and \eqref{1.29} we deduce
\be\ba{ll}\medskip
\|\bfw(t)\|\le {\rm e}^{-(\gamma-b)t}\|\bfL_1^{\beta}\bfu_0^{(1)}\|\\
\hspace*{2cm}+c_3\Int0t\Frac{{\rm e}^{-(\gamma-b)(t-s)}}{(t-s)^\beta}[(\|\bfw(s)\|+\|\bfu^{(0)}(s)\|)\|\bfw(s)\|+\|\bfw(s)\|^2]\,.\ea
\eeq{1.30}
From \lemmref{1.1} we know that the pair $(\bfw,\bfu^{(0)})$ is continuous with values in ${H}$, and so for any given $\rho>0$ there exists an interval of time $[0,\tau]$, $\tau<t^*$, such that 
\be
\sup_{t\in [0,\tau]}(\|\bfw(t)\|+\|\bfu^{(0)}(t)\|)< \rho\,,
\eeq{71}
provided 
$\|\bfw(0)\|+\|\bfu^{(0)}(0)\|<\half\rho$ (say). Thus, observing that
$$
\int_0^t\frac{{\rm e}^{-(\gamma-b)(t-s)}}{(t-s)^\beta}\,ds\le \int_0^\infty\frac{{\rm e}^{-(\gamma-b)t}}{t^\beta}\,dt:=c_0<\infty
$$
from \eqref{1.30} we get
\be
(1-c_3c_0\rho)\|{\bfw}(t)\|\le \|\bfL^{\beta}_1\bfu_0^{(1)}\|\,,\ \ t\in[0,\tau]\,.
\eeq{1.31}
We now go back to \eqref{1.20}$_2$, which, with the help of \eqref{1.27}, furnishes
$$
\|\bfu^{(0)}(t)\|\le \|\bfu^{(0)}(0)\|+c_4\int_0^t \,\big[(\|\bfu^{(1)}(s)\|+\|\bfu^{(0)}(s)\|)\|\bfu^{(1)}(s)\|+\|\hat{\bfA}^\beta\bfu^{(1)}(s)\|^2\big]\,ds\,.
$$
Thus, if we restrict to $t\in [0,\tau]$, choose $\rho=1/(2c_3c_0)$ and use \eqref{1.28}, \eqref{71} and \eqref{1.31}, the preceding inequality furnishes
\be\ba{rl}\medskip
\|\bfu^{(0)}(t)\|&\!\!\!\!\le \|\bfu^{(0)}(0)\|+c_5\Int0t \,\|\bfA^\beta\bfu(s)\|\,ds\le \|\bfu^{(0)}(0)\|+c_6\Int0t \,{\rm e}^{-b\,s}\|\bfw(s)\|\,ds \\&\!\!\!\!\le c_7 \,\|\bfA^\beta\bfu_0\|\,,\ \ t\in [0,\tau]\,.\ea 
\eeq{1.32} 
Combining \eqref{1.31} with \eqref{1.32} and  again recalling \eqref{1.28}, we  conclude in particular
\be
\|\bfA^\beta\bfu(t)\|\le c_8\,\|\bfA^\beta\bfu_0\|\,,
\eeq{1.33}
for all $t\in [0,\tau]$, namely, for as long as 
$\|\bfA^\beta\bfu(t)\|<\rho$. However, by \eqref{1.33} and a standard argument one can show that if we choose $\|\bfA^\beta\bfu_0\|$ sufficiently small, $\bfA^\beta\bfu(t)$ will never reach the boundary of the ball of radius $\rho$, implying that \eqref{1.33} must hold for all $t>0$. The proof of global existence is thus completed. Obviously, from \eqref{1.33} (valid now for all $t\ge0$), it also follows that ${\sf s}_0$ is stable, in the sense of \defref{1.1}. We shall now show the asymptotic behavior.  From what we have just proved, under the stated condition on the data and with the above choice of $\rho$, from \eqref{1.31} and \eqref{1.28} we find
\be
\|\hat{\bfA}^\beta\bfu^{(1)}(t)\|\le c_9 \,{\rm e}^{-b t} \|\hat{\bfA}^\beta\bfu^{(1)}(0)\|\,,\ \ \mbox{all $t>0$}\,.  
\eeq{1.35}
Moreover, integrating  \eqref{1.20}$_2$ between arbitrary $t_1,t_2>0$ using \eqref{1.35} and reasoning in a way similar to what we did to obtain \eqref{1.32} we infer
\be
\|\bfu^{(0)}(t_1)-\bfu^{(0)}(t_2)\|\le c_{10}\int_{t_1}^{t_2}\|\hat{\bfA}^\beta\bfu^{(1)}(s)\|ds\le c_{11}\,\|\hat{\bfA}^\beta\bfu^{(1)}(0)\|\int_{t_1}^{t_2}{\rm e}^{-b s}ds\,, 
\eeq{1.36}
from which we deduce that there exists $\bar{\bfu}^{(0)}\in \cals(\lambda)$ such that
$$
\lim_{t\to\infty}\|\bfu^{(0)}(t)-\bar{\bfu}^{(0)}\|=0\,.
$$
Plugging this information back into \eqref{1.36} with $t_2=\infty$, $t_1=t$ we infer
$$
\|\bfu^{(0)}(t)-\bar{\bfu}^{(0)}\|\le c_{12}\|\hat{\bfA}^\beta\bfu^{(1)}(0)\|\,{\rm e}^{-b t},
$$
which, once combined with \eqref{1.35}, proves the exponential rate of decay. It remains to show the instability statement. To this end we observe that, in view of  the properties of the spectrum and the estimate of the nonlinearity showed in \eqref{1.10_1}, it is easy to check that all conditions of  \cite[Theorem 5.1.3]{Henry} are satisfied, so that our statement follows from that theorem.  The proof of the proposition is thus accomplished.
\QED

\par
Combining the results of \propref{1.2} and \propref{1.1} we may conclude with the following result that furnishes necessary and sufficient conditions for the  stability  of permanent rotations.
\Bt Let $\lambda\in\{A,B,C\}$, and let ${\sf s}_0$ be the permanent rotation \eqref{1.0} with $\bfomega_0\in\cals(\lambda)-\{\0\}$. 
Then, if any of the conditions {\rm (i)--(iii)} in \propref{1.2} is satisfied, ${\sf s}_0$ is stable. Moreover, there exists $\gamma_0>0$ such that if for some $\beta\in [\frac34,1)$
$$
\|\bfA^{\beta}\bfu_0\|<\gamma_0\,,
$$
then there is $\bfu^{(0)}\equiv(\0,\bfomega^{(0)})$, $\bfomega^{(0)}\in\cals(\lambda)$, such that
all solutions to the perturbation equation \eqref{1.7} corresponding to such initial data satisfy
$$
\|\bfA^{\beta}\bfu(t)-\bfu^{(0)}\|\le C\,\|\bfA^{\beta}\bfu_0\|\,{\rm e}^{-\kappa\,t}\,,\ \ \mbox{all $t\ge 0$}\,,
$$ 
for some constants $C,\kappa>0$. Conversely, if any of the conditions {\rm (iv)--(vi)} in \propref{1.2} holds, the permanent rotation ${\sf s}_0$ is unstable. 
\ET{1.1}
\Br The result of \theoref{1.1} provides, in particular, a sharp and rigorous explanation of Lord Kelvin's experiment, mentioned in the introductory section. Actually, this experiment shows that rotations occurring around the shorter axis (=\,maximum moment of inertia) of a prolate spheroid filled with water  are stable. while those around the other central axes are unstable. As a matter of fact, our result implies a much more general phenomenon, namely, that {\em whatever the  shape of the body and cavity}, the {\em only stable rotations} are those occurring along the axis with respect to which the moment of inertia is a maximum, all others being unstable.  
\ER{3.2}
\renewcommand{\theequation}{{4}.\arabic{equation}}
\setcounter{equation}{0}
\setcounter{section}{4}\setcounter{lemm}{0}
\setcounter{theo}{0}
\section*{\normalsize 4. Asymptotic Behavior for Large Data. A Full Proof of Zhukovsky Conjecture.}
In view of \theoref{1.1} we may deduce that under any of the assumptions (i)--(ii) stated in \propref{1.2} all solutions to \eqref{1.1}--\eqref{1.2} that belong to a suitable function class and with initial data ``sufficiently close" to a permanent rotation  \eqref{1.3} in the eigenspace $\cals(\lambda)$ must converge exponentially fast to (another) permanent rotation that lies in the same eigenspace.  Objective of this section is to show that, in fact, the same conclusion holds in the more general  class of {\em weak} solutions to \eqref{1.1}--\eqref{1.2} and for data that not only are less regular, but also not necessarily ``close" to a permanent rotation. 
These results provide, in particular, a positive answer to an important question that was left open in \cite{DGMZ}, thus providing a completely affirmative answer to the  {\em Zhukovsky conjecture} mentioned in the introduction.  
\setcounter{defn}{0}
\vspace*{.1mm}\par
We begin to recall the definition of weak solution.
\Bd A pair $(\bfv,\bfomega)$ is a {\em weak solution} to \eqref{1.1}--\eqref{1.2} if it satisfies the following properties
\begin{itemize}
  \item[{\rm (i)}] $\bfv\in C_w([0,\infty);L^2_\sigma(\mathcal C))\cap L^\infty(0,\infty;L^2_\sigma(\mathcal C))\cap L^2(0,\infty;W_0^{1,2}(\mathcal C))$\,; 
  \item[{\rm (ii)}] $\bfomega\in C([0,\infty))$\,,\ $(\bfomega-\bfa)\in W^{1,\infty}((0,\infty))$\,; 
  \item [{\rm (iii)}] Strong energy inequality:  
$$  
\mathcal E(t)+\nu\int_s^t\|\nabla\bfv(\tau)\|_2^2\le \mathcal E(s)
$$
for $s=0$, for a.a. $s>0$ and all $t\ge s$,  where 
$$
\mathcal E:=\half\left(\|\bfv\|^2_2+\bfomega\cdot\mathbb I\cdot\bfomega-2\bfa\cdot\mathbb I\cdot\bfomega\right)\,,
$$
and, we recall,  $\bfa$ is defined in \eqref{1.5};
\item[{\rm (iv)}] $(\bfv,\bfomega)$ obey \eqref{1.1}--\eqref{1.2} in the sense of distribution\,. 
\end{itemize}
\EED{4.1}
\Br
Employing the important property that $\mathcal E$ is positive definite in the variables $(\bfv,\bfomega)$, one can show  \cite[\S 3.2]{M} that the class, $\mathscr C_w$, of weak solutions is not empty. Moreover, any $(\bfv,\bfomega)\in\mathscr C_w$ is unique in the class of those $(\bfw,\varpi)\in\mathscr C_w$ such that $\bfw\in L^q(0,\tau;L^r(\calc))$, $3r^{-1}+2q^{-1}=1$ \cite[\S 3.4]{M}. 
\ER{4.3}

We need the following preliminary result.

\Bl Let $(\bfv,\bfomega)$ be a weak solution to the problem \eqref{1.1}--\eqref{1.2} corresponding to data $(\bfv_0,\bfomega_0)\in L^2_\sigma(\calc)\times \real^3$. Then, 
there exists $t_0>0$ such that for all $t>t_0$ the solution becomes regular and, in particular,  satisfies the following conditions,
\be\ba{ll}\medskip  
\bfv\in W^{1,\infty}(t_0,\infty; L^2_\sigma(\calc))\cap C([t_0,\infty);W^{1,2}_0(\calc))\cap L^\infty([t_0,\infty);W^{2,2}(\calc))\,,\\
\bfomega\in C^1([t_0,\infty);\real^3)\,.
\ea
\eeq{2.0_0}
Moreover,
\be
\lim_{t\to\infty}\|\bfA_0^\alpha\bfv(t)\|_2=0\,,\ \ \mbox{for all $\alpha\in [0,1)$}\,.
\eeq{2.0}
\EL{2.1}
{\em Proof.} Setting
\be
\bfomega_\infty:=\bfomega-\bfa\,,
\eeq{schu1}
we have that the governing equations \eqref{1.1}--\eqref{1.2} can be rewritten as follows
\be\ba{cc}\medskip
\left.\ba{ll}\medskip
\bfv_t+\bfv\cdot\nabla\bfv+(\dot{\bfomega}_\infty+\dot{\bfa})\times\bfx+2(\bfomega_\infty+\bfa)\times\bfv=\nu\Delta\bfv-\nabla p\\
\Div\bfv=0\ea\right\}\ \ \mbox{in $\calc\times (0,\infty)$}\\
\bfv(x,t)=\0\ \ \mbox{at $\partial\Omega$}
\ea\eeq{2.1}
and
\be
\mathbb I\cdot\dot{\bfomega}_\infty+(\bfomega_\infty+\bfa)\times\mathbb I\cdot\bfomega_\infty=\0\,.
\eeq{2.2}
From \cite[Proposition 1]{DGMZ} we know already that, for any given  weak solution $(\bfv,\bfomega_\infty)$ to \eqref{2.1}--\eqref{2.2} there is $t_0>0$ such that for all $t\ge t_0$ the solution becomes strong. Precisely, for all $T>0$:
\be\ba{ll}\medskip
\bfv\in C([t_0,\infty);W_0^{1,2}(\calc))\cap L^2(t_0,t_0+T;W^{2,2}(\calc))\,,\\ \bfv_t\in L^2(t_0,t_0+T;L^{2}(\calc))\,,\ \ \bfomega_\infty\in C^1([t_0, \infty); \real^3)\,,
\ea\eeq{2.3}
and there exists $p\in L^2(t_0,t_0+T;W^{1,2}(\calc))$,  such that $(\bfv,p,\bfomega_\infty)$ 
satisfy \eqref{2.1}--\eqref{2.2}, a.e. in $\mathcal C\times (0,\infty)$.
Finally, 
\be
	\lim_{t\to\infty} \Vert \bfv (t)\Vert_{1,2} = 0. 
	\eeq{2.4}
In view  of \eqref{2.3}$_1$, and \eqref{2.4} we deduce that for any $\eta>0$ there exists some time $t'\ge t_0$ such that  the $W^{2,2}$-norm of $\bfv$ is finite at $t'$, while the $W^{1,2}$-norm is less than $\eta$ for all $t\ge t'$. Without loss of generality, we may take $t'=t_0$ and shall thus suppose 
\be
\|\bfv(t_0)\|_{2,2}<\infty\,,\ \ \|\bfv(t)\|_{1,2}<\eta\,, \ \,\mbox{for all $t\ge t_0$}.
\eeq{2.5}
Our next goal is to construct a solution to \eqref{2.1}--\eqref{2.2} for $t\ge t_0$,  corresponding to the initial  data $(\bfv(t_0),\bfomega_\infty(t_0))$. It will be achieved  by an argument  analogous to that presented in \cite[Theorem 4.1]{GaMa}. To this end, we will prove some basic {\em a priori} estimates for solutions to \eqref{2.1}--\eqref{2.2}. Taking the time-derivative of both sides of \eqref{2.1}, then dot-multiplying the resulting equation by $\bfv_t$ and integrating by parts over $\calc$, we get (formally)
\be
\frac12\ode{}t\left(\|\bfv_t\|_2^2-\dot{\bfa}\cdot\mathbb I\cdot\dot{\bfa}\right)+\nu\|\nabla\bfv_t\|_2^2=-(\bfv_t\cdot\nabla\bfv,\bfv_t)+\ddot{\bfomega}_\infty\cdot\mathbb I\cdot \dot{\bfa} - 2((\dot{\bfomega}_\infty+\dot{\bfa})\times\bfv,\bfv_t)\,.
\eeq{2.6}
Moreover, from \cite{DGMZ} we know that $(\bfv,\bfomega_\infty)$ satisfy the energy balance equation (energy equality)
\be
\frac12\ode{}t\left(\|\bfv\|_2^2-{\bfa}\cdot\mathbb I\cdot{\bfa}+\bfomega_\infty\cdot\mathbb I\cdot\bfomega_\infty\right)+\nu\|\nabla\bfv\|_2^2=0\,.
\eeq{2.7}
From \eqref{2.2} it follows that
\be
|\dot{\bfomega}_\infty|\le c_1\,(|\bfomega_\infty|^2+\|\bfv\|_2^2)\,,
\eeq{s1}
so that
\be
|2((\dot{\bfomega}_\infty+\dot{\bfa})\times\bfv,\bfv_t)|\le c_1\big[(|\bfomega_\infty|^2+\|\bfv\|_2^2)\|\bfv_t\|_2+\|\bfv\|_2\|\bfv_t\|_2^2\big]\,.
\eeq{s2}
Furthermore, taking the time-derivative of both sides of \eqref{2.2} and dot-multiplying the resulting equation by $\dot{\bfa}$ produces 
$$
\ddot{\bfomega}_\infty\cdot\mathbb I\cdot \dot{\bfa}=-\dot{\bfomega}_\infty\times\mathbb I\cdot\bfomega_\infty\cdot\dot{\bfa}-(\bfomega_\infty+\bfa)\times\mathbb I\cdot\dot{\bfomega}_\infty\cdot\dot{\bfa}\,.
$$
Therefore, employing \eqref{s1} we show
\be
|\ddot{\bfomega}_\infty\cdot\mathbb I\cdot \dot{\bfa}|\le c_2\,\big[(|\bfomega_\infty|^2+\|\bfv\|_2^2)(|\bfomega_\infty|+\|\bfv\|_2)\|\bfv_t\|_2\,.
\eeq{s3}
Now, from
  \eqref{2.7} we obtain rather simply
\be
\|\bfv(t)\|_2+|\bfomega_\infty(t)|\le D\,,
\eeq{2.8}
where $D$, here and in the following, denotes a constant depending only on the initial data. Using \eqref{2.8} and \eqref{s2} we then deduce
\be
|2((\dot{\bfomega}_\infty+\dot{\bfa})\times\bfv,\bfv_t)|\le c_3(1+\|\bfv\|_2\|\bfv_t\|_2^2)
\eeq{2.10}
with $c_3$ depending on $D$. Furthermore,   employing in  \eqref{s3} the bound  \eqref{2.8} and the Cauchy-Schwarz inequality, we show for arbitrary $\varepsilon>0$
\be
|\ddot{\bfomega}_\infty\cdot\mathbb I\cdot \dot{\bfa}|\le c_4+\varepsilon\|\bfv_t\|_2^2
\eeq{2.11}
with $c_4$ depending on $D, \varepsilon$. Furthermore, by Cauchy-Schwarz  inequality and Sobolev embedding theorem, it follows that
\be
|(\bfv_t\cdot\nabla\bfv,\bfv_t)|\le \|\bfv_t\|_4^2\|\nabla\bfv\|_2\le \|\nabla\bfv_t\|_2^{\frac32}\|\bfv_t\|_2^{\frac12}\|\nabla\bfv\|_2\le \varepsilon\,\|\nabla\bfv_t\|_2^2+c_5\,\|\bfv_t\|_2^2\|\nabla\bfv\|_2^4\,.
\eeq{2.12}
whereas, by \eqref{2.7}, \eqref{s1} and \eqref{2.8}, we also show
\be
\|\nabla\bfv\|_2^2\le c_4 (\|\bfv_t\|_2+1)\,.
\eeq{2.13}
As a result, combining \eqref{2.6}, \eqref{2.10}--\eqref{2.13}, and by choosing $\varepsilon=\nu/2$ in \eqref{2.12}, we infer
\be
\ode{}t\left(\|\bfv_t\|_2^2-\dot{\bfa}\cdot\mathbb I\cdot\dot{\bfa}\right)+\nu\|\nabla\bfv_t\|_2^2\le c_5(1+\|\bfv_t\|_2^4)\,.
\eeq{2.14}
If we now  set
$$
E_1:=\frac12(\|\bfv_t\|_2^2-\dot{\bfa}\cdot\mathbb I\cdot\dot{\bfa}) 
$$
from \eqref{1.44} we deduce
\be
c_0\|\bfv_t\|_2^2\le 2 E_1\le \|\bfv_t\|_2^2. 
\eeq{2.14_1}
Thus \eqref{2.14}  furnishes, in particular,
$$
\ode{}tE_1\le c_6\,(1+E_1)^4\,.
$$
Integrating this inequality and taking into account \eqref{2.14_1}, we easily show that there exists 
\be
\tau\ge c_7/(\|\bfv_t(t_0)\|_2+1)^3\,,
\eeq{2.15}  
and continuous functions $H_i$, $i=1,2$, in $[t_0,t_0+\tau)$ such that 
\be
\|\bfv_t(t)\|_2\le H_1(t)\,,\ \ \int_{t_0}^t\|\nabla\bfv_s(s)\|_2^2\le H_2(t)\ \ t\in[t_0,t_0+\tau)\,.
\eeq{2.16}
Using \eqref{2.1}$_1$ and proceeding exactly as in \cite[Theorem 4.1]{GaMa}
we can prove that
\be
\|\bfv_t(t_0)\|_2\le c_8\,[{\sf P}(\|\bfv(t_0)\|_{2,2})+1] 
\eeq{2.16_1}
where ${\sf P}(r)$ is a polynomial in $r$ with ${\sf P}(0)=0$. Thus, in view of \eqref{2.5} and \eqref{2.15}, this implies, on the one hand, that $\|\bfv_t(t_0)\|_2$ is well defined and, on the other hand the following estimate on $\tau$
\be
\tau\ge {c}{/ ({\sf P}(\|\bfv(t_0)\|_{2,2})+1)^3}.
\eeq{chip}
Estimates \eqref{2.16} along with \eqref{2.13} allow us to establish by the classical Faedo-Galerkin method the existence of a solution $(\tilde{\bfv},\tilde{\bfomega})$ to \eqref{2.1}--\eqref{2.2} in the time interval $[t_0,t_0+\tau)$ with the following properties valid for all $t_1\in(0,\tau)$
$$
\ba{ll}\medskip
\tilde{\bfv}\in W^{1,\infty}([t_0,t_1]);L^2(\calc))\cap L^\infty(t_0,t_1;W^{1,2}(\calc))\,,\\ \tilde{\bfv}_t\in L^2(t_0,t_1;W^{1,2}(\calc))\,,\ \ \tilde{\bfomega}_\infty\in C^1([t_0, t_1]); \real^3)\,;
\ea
$$
for details, see \cite[Chapter 4]{M}. 
By the uniqueness result recalled in \remref{4.3}, we must have $(\tilde{\bfv},\tilde{\bfomega})=({\bfv},{\bfomega})$ in $[t_0,t_0+\tau)$.
Using \eqref{chip} and arguing in a way entirely analogous to \cite[Theorem 4.1]{GaMa}, one can then show that either $\tau=\infty$, or else  $\|\bfv(t)\|_{2,2}$ becomes unbounded in a (left) neighborhood of $\tau$. In view of\eqref{EmNa_0},
the above is equivalent to say that $\|\bfA_0\bfv(t)\|_2$ becomes unbounded. 
We shall prove that 
\be
\sup_{t\in [t_0,t_0+\tau)}\|\bfA_0\bfv(t)\|_{2,2}\le M<\infty\,,
\eeq{2.17}
and thus conclude $\tau=\infty$.
To secure \eqref{2.17}, we go back to \eqref{2.6} and employ the estimates \eqref{2.10}--\eqref{2.12} on its right-hand side. We thus obtain
$$
\ode{}tE_1+\nu\|\nabla\bfv_t\|_2^2\le c_9+\varepsilon\,\|\bfv_t\|_{1,2}^2 +c_{10}\,(\|\bfv\|_2^2+\|\nabla\bfv\|_2^4)\|\bfv_t\|_2^2
$$
Observing that, by Poincar\'e inequality and \eqref{2.14_1}, 
\be
\|\nabla\bfv\|_2^2\ge  c_{11}\,\|\bfv_t\|_2^2\ge 2c_{11}E_1\,,
\eeq{2.18}
we may take $\varepsilon$ sufficiently small to obtain 
$$
\ode{}tE_1+c_{12}E_1\le c_9 +2c_{10}\,(\|\bfv\|_2^2+\|\nabla\bfv\|_2^4)E_1
$$
However, by \eqref{2.5}, $\eta$ can be chosen in such a way that
$$
2c_{10}\,(\|\bfv\|_2^2+\|\nabla\bfv\|_2^4)<\half c_{12}\,,
$$
so that the preceding inequality becomes
$$
\ode{}tE_1+c_{13}E_1\le c_9
$$
Integrating this differential equation between $t_0$ and $t<\tau$, and taking into account \eqref{2.16_1}, \eqref{2.5} we show
\be
\|\bfv_t(t)\|_2\le D\,,\ \ t\in [t_0,\tau)\,.
\eeq{2.22}
Next, dot-multiplying both sides of \eqref{2.1} by ${\rm P}\Delta$ and using Schwarz inequality, we also show
\be
\|\bfA_0\bfv(t)\|_2\le c_{14}\left[\|\bfv_t(t)\|_2+(|\bfomega_\infty(t)|+\|\bfv(t)\|_2)\|\bfv(t)\|_2+|\dot{\bfomega}_\infty(t)|+\|\bfv(t)\cdot\nabla\bfv(t)\|_2\right]\,.
\eeq{2.23}
The following inequality is well known (e.g., \cite[eq. (3.22)]{GaMa})
$$
\|\bfv\cdot\nabla\bfv\|_2\le c_{14}\|\nabla\bfv\|_2^3+\varepsilon\,\|\bfA_0\bfv\|_{2}
$$
with arbitrary $\varepsilon>0$ and $c_{14}\to\infty$ as $\varepsilon\to 0$.
Thus, using the latter with $\varepsilon=1/(2c_{14})$ into \eqref{2.23} delivers
\be
\|\bfA_0\bfv(t)\|_2\le 2c_{14}\left[\|\bfv_t(t)\|_2+(|\bfomega_\infty(t)|+\|\bfv(t)\|_2)\|\bfv(t)\|_2+|\dot{\bfomega}_\infty(t)|+\|\nabla\bfv(t)\|_2^3\right]\,,
\eeq{GaMa}
which, by virtue of \eqref{2.5}, \eqref{s1},  \eqref{2.8},  and \eqref{2.22} allows us to show \eqref{2.17}. As a result, \eqref{2.17} holds with $\tau=\infty$, namely
\be
\sup_{t\ge t_0}\|\bfA_0\bfv(t)\|_2\le M\,.
\eeq{2.25} 
This property, in conjunction with \eqref{2.3} and \eqref{2.22}, shows that $(\bfv,\bfomega)$ is in the class \eqref{2.0_0}.
Finally, we recall the well-known inequality (e.g., \cite[p. 28]{Henry})
$$
\|\bfA_0^\alpha\bfv(t)\|_2\le \|\bfA_0\bfv(t)\|_2^\alpha\|\bfv(t)\|_2^{1-\alpha}\,,\ \ \alpha\in [0,1]
$$
which, with the help of \eqref{2.4} and \eqref{2.25} implies \eqref{2.0}. The proof is completed.
\QED
\Br We wish to observe that, by using arguments similar to those employed in the proof of the previous lemma, one can show that the weak solution  possesses  regularity properties,  in the time variable, even stronger than that stated in \eqref{2.0_0}$_1$. More precisely, one can show $\bfv\in W^{k,2}((t_0,T); W^{2,2}(\calc))$, for all $k \ge 1$ and $T>0$ (e.g., \cite[Lemma 5.5]{IVP}). 
\ER{0.1}

We are now in a position to prove the  main result of this section that gives a complete and positive answer to {\em Zhukovsky's conjecture}. 
\Bt Let ${\mathfrak s}:=(\bfv,\bfomega)$ be a weak solution to \eqref{1.1}--\eqref{1.2} corresponding to initial data $(\bfv(0),\bfomega(0))\in L^2_\sigma(\calc)\times\real^3$. Then, there exists $t_0=t_0(\mathfrak s)$ such that, for all $t\ge t_0$, $\mathfrak s$ becomes smooth and, precisely,  lies in the function class defined by \eqref{2.1}. Moreover,
\be
\lim_{t\to\infty}\|\bfA_0^{\alpha}\bfv(t)\|_2=0,
\eeq{2.28}
for all $\alpha\in [0,1)$, and there exists $\bar{\bfomega}\in\real^3$ such that 
\be
\lim_{t\to\infty}|\bfomega(t)-\bar{\bfomega}|=0\,,
\eeq{2.29}
the rate of decay in both \eqref{2.28} and \eqref{2.29} being exponential.
\ET{2.1}
{\em Proof.} The first part of the theorem that includes \eqref{2.28}, has been proved in \lemmref{2.1}. Concerning the property related to \eqref{2.29}, in \cite[Theorem 4]{DGMZ} it is shown that it is always true with $\bar{\bfomega}=\0$ in the (physically irrelevant) case that the initial total angular momentum $\bfM(0)$ vanishes (see \eqref{1.2}). If, however, $\bfM(0)\neq\0$ (which necessarily implies $\bar{\bfomega}\neq\0)$, in \cite[loc.cit.]{DGMZ} the above property is secured only  when either $A\le B<C$ or  $A=B=C$. Thus, the case $A<B=C$ is left open. In fact, in such circumstance it is only shown that either, for some $\bar{\omega}\in \real$,
$$ 
\lim_{t\to\infty}|\bfomega(t)-\bar{\omega}\bfe_1|=0\,,
$$
in which case \eqref{2.29} is proved,  or else
\be
\lim_{t\to\infty}{\rm dist}(\bfomega(t),\cals(\lambda))=0\,,
\eeq{2.30}
with $\lambda\equiv B=C$. Therefore, in the occurrence of \eqref{2.30}, we cannot assert the validity of \eqref{2.29}. However,  with the help of \theoref{1.1} we will now show that also in the case \eqref{2.30} there must be $\bar{\bfomega}$ ($\neq\0$) for which \eqref{2.29} holds. Actually, from \eqref{2.30} it follows that we can find $\bfomega_0\in\cals(\lambda)$, $\bfomega_0\neq\0$, and an unbounded sequence of times $\{t_n\}$ such that
$$
\lim_{n\to\infty}|\bfomega(t_n)-\bfomega_0|=0\,.
$$
In view of the latter and \eqref{2.28}, there is  $\tau_0\ge t_0$ such that 
$$
\|\bfA_0^{\beta}\bfv(\tau_0)\|_2+|\bfomega(\tau_0)-\bfomega_0|<\gamma_0\,,
$$
where $\beta$ and $\gamma_0$ are the constants introduced in \theoref{1.1}. Because of the uniqueness property (see \remref{4.3}), our weak solution will then coincide with the solution of that theorem. In particular, since the assumption $A<B=C$ is exactly case (ii) of \propref{1.2},  from \theoref{1.1} we conclude the existence of $\bar{\bfomega}\in\cals (\lambda)$ for which \eqref{2.29} holds. The exponential decay rate, in all cases, is ensured by \eqref{2.28}, \eqref{2.29} and \theoref{1.1}. The proof is completed.

\QED

\renewcommand{\theequation}{{5}.\arabic{equation}}
\setcounter{equation}{0}
\setcounter{section}{5}\setcounter{lemm}{0}
\setcounter{theo}{0}
\section*{\normalsize 5.  Attainability of Permanent Rotations} 
One of the significant questions posed by \theoref{2.1} concerns around which central axis the terminal permanent rotation will occur. Notice that this problem is further reinforced by the fact that, since we are dealing with {\em weak} solutions for which uniqueness is {\em not} known, it may happen, in principle, that two solutions corresponding to the {\em same} initial data may converge, eventually, to permanent rotations occurring around two {\em different} central axes. 

The above question has been analyzed in some detail in \cite[Section 8]{DGMZ}, where it is shown  for an open set of ``large" data that, provided $A\le B<C$, the terminal permanent rotation, ${\sf r}_0$, will take place around the $\bfe_3$-axis, namely, the axis with maximum moment of inertia. However, this result leaves open two important aspects. In the first place, it does not allow us to draw an analogous conclusion if $A<B=C$, that is, ${\sf r}_0$ will occur around an axis spanned by $\{\bfe_2,\bfe_3\}$, as somehow expected on the basis of \theoref{2.1}. Moreover, it does not provide any rate of decay. In this regard, both lab  \cite{WJ} and numerical \cite[Section 9.1]{DGMZ} experiments show that, after a transient interval of time, whose length depends (inversely) on the kinematic viscosity coefficient, the motion of the coupled system almost abruptly converges to a permanent rotation. \par Objective of this section is to fill these two gaps. In particular, we shall show that when the (relative) velocity of the liquid becomes ``sufficiently small" and the angular velocity of the body is ``sufficiently close" to that of the corresponding terminal permanent rotation, both quantities must decay exponentially fast to their respective limits, which thus also explains the ``abrupt" convergence mentioned before. 
\Bt Let $(\bfv,\bfomega)$ be a weak solution in the sense of \defref{1.1}, corresponding to data $(\bfv_0,\bfomega_0)$, and let $E=E(t)$ the functional defined in \eqref{ene}.
The following properties hold.\footnote{We assume  $\bfM(0)\neq\0$, otherwise the motion of the coupled system is physically irrelevant; see \remref{3.1}.}
\begin{itemize}
\item[{\rm (a)}] If $A=B<C$, assume
$$
E(0)\le\frac{(C-A)\,C}{2 A}\omega_{03}^2(0)\,,
$$
whereas, if $A<B<C$, assume  
$$\ba{ll}\medskip
E(0)+\Frac{A}{2B}(B-A)\omega_{01}^2(0)\le \Frac C{2B}(C-B)\omega_{03}^2(0)\,,\\
E(0)\le\Frac{B}{2A}(B-A)\omega_{02}^2(0)+\Frac C{2A}(C-A)\omega_{03}^2(0)\,.
\ea
$$
Then, there exists $\bar{\bfomega}\in \cals(C)-\{\0\}$ such that
\be\lim_{t\to\infty}\big(\|\bfA_0^\alpha\bfv(t)\|_2
+|\bfomega(t)-\bar{\bfomega}|\big)=0\,,\ \ \mbox{all $\alpha\in [0,1)$\,.}
\eeq{5.1} 
\item[{\rm (b)}] 
If $A<B=C$, assume 
\be
E(0)\le \frac{B(B-A)}{2 A}(\omega_{02}^2(0)+\omega_{03}^2(0))\,.
\eeq{5.2}
Then, \eqref{5.1} holds for some $\bar{\bfomega}\in\cals (B)\ (\equiv\cals(C))-\{\0\}$.
\item[{\rm (c)}] Under the assumptions stated in {\rm (a)} and {\rm (b)}, there exist positive constants $t_0,C$ and $\gamma$, depending at most on $\nu,\calc, A,B,C$ and $(\bfv_0,\bfomega_0)$, such that 
$$
\|\bfv_t(t)\|_2+\|\bfv(t)\|_{2,2}+|\dot{\bfomega}(t)|+|\bfomega(t)-\bar{\bfomega}|\le C\,{\rm e}^{-\gamma t}\,,\ \ \mbox{for all $t\ge t_0$.}
$$
\end{itemize}
\ET{5.1}
{\em Proof.} We begin to prove the statement in (a). Under the given assumptions on the initial data, in \cite[Theorem 5]{DGMZ} it is shown that
$$
\lim_{t\to\infty}|\bfomega(t)-\bar{\bfomega}|=0\,.
$$
The claim then follows from \lemmref{2.1}. In  \cite[Theorem 5]{DGMZ} it has also proved that if \eqref{5.2} holds, necessarily 
$$
\lim_{t\to\infty}|\omega_1(t)|=0\,.
$$
As a result, by \theoref{2.1} we deduce that $\bfomega(t)$ must converge, as $t\to\infty$, to some $\bar{\bfomega}\in\cals(B)\, (\equiv\cals(C))$, which proves the statement in (b). We shall next show (c). From \theoref{1.1} we already know, in particular, that
\be 
\|\bfv(t)\|_{1,2}+|\bfomega(t)-\bar{\bfomega}|\le C\, {\rm e}^{-\kappa t}\,,\ \ \mbox{for all $t\ge t_0$}\,,
\eeq{5.3}
where, here and in what follows, $C$ denotes a generic constant whose value may change from line to line, that depends, at most, on  $\nu,\calc$, the central moments of inertia of the coupled system body-liquid and the initial data. Observing that $\bar{\bfomega}\times\mathbb I\cdot\bar{\bfomega}=\0$, from \eqref{2.2} we have 
\be
\mathbb I\cdot\dot{\bfomega}_\infty=-(\bfw+\bfa)\times\mathbb I\cdot\bfomega_\infty-\bar{\bfomega}\times\mathbb I\cdot\bfw\,,
\eeq{5.4}
with $\bfw:=\bfomega_\infty-\bar{\bfomega}$. Therefore, from \eqref{2.8}, \eqref{5.3}, and \eqref{2.8} we conclude
\be
|\dot{\bfomega}_\infty(t)|\le C\,{\rm e}^{-\kappa t}\,,\ \ t\ge t_0\,.
\eeq{5.5} 
By Schwarz inequality, \eqref{5.3} and \eqref{5.5}, we show
\be\ba{rl}\medskip
|2((\dot{\bfomega}_\infty+\dot{\bfa})\times\bfv,\bfv_t)|\le &\!\!\!\!2|\dot{\bfomega}_\infty|\,\|\bfv\|_2\|\bfv_t\|_2+C\|\bfv\|_2\|\bfv_t\|_2^2\le |\dot{\bfomega}_\infty|^2\|\bfv\|_2+ C\,\|\bfv\|_2\|\bfv_t\|_2^2\\
\le &\!\!\!\! C\,{\rm e}^{-\kappa t}(1+\|\bfv_t\|_2^2)\,.
\ea
\eeq{5.6}
Moreover, taking the time-derivative of both sides of \eqref{5.4} and dot-multiplying the resulting equation by $\dot{\bfa}$, we get
$$
\ddot{\bfomega}_\infty\cdot\mathbb I\cdot \dot{\bfa}=-\dot{\bfw}\times\mathbb I\cdot\bfomega_\infty\cdot\dot{\bfa}-(\bfw+\bfa)\times\mathbb I\cdot\dot{\bfomega}_\infty\cdot\dot{\bfa}-\bar{\bfomega}\times \mathbb I\cdot\dot{\bfomega}_\infty\cdot\dot{\bfa}\,.
$$
Thus, employing in this relation Cauchy-Schwarz inequality along with \eqref{5.3} and \eqref{5.5} we obtain
\be
|\ddot{\bfomega}_\infty\cdot\mathbb I\cdot \dot{\bfa}|\le C\,{\rm e}^{-\kappa t}(1+\|\bfv_t\|_2^2)\,.
\eeq{5.7}
We next observe that, in view of \theoref{1.1} and \lemmref{2.1}, the generic weak solution becomes smooth for sufficiently large $t$ (see also \remref{0.1}). Therefore, in particular,  
 \eqref{2.6} holds for all such instant of times. Now, the latter in conjunction with \eqref{2.12} (with $\varepsilon=\nu/2$) , \eqref{5.3}, \eqref{5.6} and \eqref{5.7}, entails that the weak solution satisfies, for all large $t$, the following inequality
$$
\ode{}t \big(\|\bfv_t\|^2_2-\dot{\bfa}\cdot\mathbb I\cdot\dot{\bfa}\big)+\nu\,\|\nabla\bfv_t\|_2^2\le C\,{\rm e}^{-\kappa t}(1+\|\bfv_t\|_2^2)\,.
$$
By Poincar\'e inequality and \eqref{2.14_1}, this relation allows us to conclude that for all $t>\kappa^{-1}\ln(2C/\nu):=\tau$, we must have
\be
\ode{}tE_1+\gamma\,E_1\le C \,{\rm e}^{-\kappa t}\,,
\eeq{5.8}
for suitable constant $\gamma=\gamma(\nu,\Omega)>0$. Notice that, without loss of generality, we can  assume  $\gamma<\kappa$, because if \eqref{5.8} holds for some $\gamma$, it continues to hold for $\gamma_1<\gamma$.   If we integrate \eqref{5.8} from $\tau$ to $t>\tau$, we show 
$$
E_1(t)\le E_1(\tau)\,{\rm e}^{-\gamma t}+C\,{\rm e}^{-\gamma t}\int_\tau^t{\rm e}^{-(\kappa-\gamma) s}ds\,, \ \ t\ge\tau
$$
from which, with the help of  \eqref{2.14_1} and observing that $\kappa>\gamma$, we deduce
\be
\|\bfv_t(t)\|_2^2\le C\, {\rm e}^{-\gamma t}
\eeq{5.9}
If we now consider inequality \eqref{GaMa} and estimate its right-hand side with the help of \eqref{5.3}, \eqref{5.5} and \eqref{5.9}. we arrive at
$$
\|\bfA_0\bfv(t)\|_2\le C\, {\rm e}^{-\gamma t}\,,
$$
which by \eqref{EmNa_0} allows us to conclude
$$
\|\bfv(t)\|_{2,2}\le C\, {\rm e}^{-\gamma t}\,.
$$
The statement in part (c) then follows from the latter, \eqref{5.4}, \eqref{5.5}, \eqref{5.9}, and the theorem is completely proved.
\QED

\medskip\par\noindent{\bf Acknowledgments.} Work partially supported by NSF grant DMS--1614011. I would like to thank Professor Jan Pr\"uss and Mr. Jan A. Wein for inspiring conversations.

\ed

\ed

and, moreover, by Poincar\'e inequality, 
\be
\int_0^\infty\|\bfv(t)\|_2^2dt\le M\,.
\eeq{1.49}
Now, \eqref{1.48} 
along with \eqref{1.43}$_4$, furnishes  
$$
|\dot{\bfomega}_*(t)|\le M_1\ \ \mbox{ all $t>0$,}
$$
for another constant $M_1>0$. Replacing this information on the right-hand side of \eqref{1.44} and using \eqref{1.45} and Schwarz and Poincar\'e inequalities, we show
$$
\ode{E}t+c_1\, E\le c_2\,E^{\frac12}\,,
$$
which, by a generalized form of Gronwall's lemma \cite[Lemma 2.1]{GMZ} and \eqref{1.45}, \eqref{1.49} implies
\be
\lim_{t\to\infty}\|\bfv(t)\|_2=0\,.
\eeq{1.50}
From \eqref{1.43}$_5$, \eqref{1.48} and \eqref{1.50} we deduce that the orbits generated by the solutions to \eqref{1.37} through any initial data $\bfu_0$ are compact and, therefore, the $\omega$-limit set is not empty and, in particular, invariant. By \eqref{1.5}, $\bfv\equiv \0$ on this set so that, from \eqref{1.43}$_1$ we derive $\boms=\bar{\bfomega}=\textrm{const.}$  and, from \eqref{1.43}$_4$ that $\bar{\bfomega}$ must satisfy
$$
\bfomega_0\times\mathbb I\cdot\bar{\bfomega}+\bar{\bfomega}\times\mathbb I\cdot\bfomega_0=\0\,,
$$
implying, by \lemmref{1.0}, $\bar{\bfomega}\in \cals(\lambda)$ which, since $\bar{\bfomega}\in \cals(\lambda)^\perp$, in turn entails $\bar{\bfomega}=\0$. Thus, by definition of $\omega$-limit set, we conclude
$$
\lim_{t\to\infty}|\boms(t)|=0
$$
that, once combined with \eqref{1.50}, proves \eqref{1.38}. The proof of the stability property is thus completed.


If, however, $\bfM(0)\neq\0$ (which necessarily implies $\bar{\bfomega}\neq\0)$, in \cite[loc.cit.]{DGMZ} the above property is secured only  when either $A\le B<C$ or  $A=B=C$. Thus, the case $A<B=C$ is left open. In fact, in such circumstance it is only shown that either, for some $\bar{\omega}\in \real$,
$$ 
\lim_{t\to\infty}|\bfomega(t)-\bar{\omega}\bfe_1|=0\,,
$$
in which case \eqref{2.29} is proved,  or else
\be
\lim_{t\to\infty}{\rm dist}(\bfomega(t),\cals(\lambda))=0\,,
\eeq{2.30}
with $\lambda\equiv B=C$. Therefore, in the occurrence of \eqref{2.30}, we cannot assert the validity of \eqref{2.29}. However,  with the help of \theoref{1.1} we will now show that also in the case \eqref{2.30} there must be $\bar{\bfomega}$ ($\neq\0$) for which \eqref{2.29} holds. Actually, from \eqref{2.30} it follows that we can find $\bfomega_0\in\cals(\lambda)$, $\bfomega_0\neq\0$, and an unbounded sequence of times $\{t_n\}$ such that
$$
\lim_{n\to\infty}|\bfomega(t_n)-\bfomega_0|=0\,.
$$
In view of the latter and \eqref{2.28}, there is  $\tau_0\ge t_0$ such that 
$$
\|\bfA_0^{\beta}\bfv(\tau_0)\|_2+|\bfomega(\tau_0)-\bfomega_0|<\gamma_0\,,
$$
where $\beta$ and $\gamma_0$ are the constants introduced in \theoref{1.1}. Because of the uniqueness property (see \remref{4.3}), our weak solution will then coincide with the solution of that theorem. In particular, since the assumption $A<B=C$ is exactly case (ii) of \propref{1.2},  from \theoref{1.1} we conclude the existence of $\bar{\bfomega}\in\cals (\lambda)$ for which \eqref{2.29} holds. The theorem is completely proved.

\